\documentclass[12pt,reqno]{amsart}
\usepackage{amssymb}
\usepackage{amsmath}
\usepackage{bm}
\usepackage{verbatim}
\makeatletter
\@addtoreset{equation}{section}
\makeatother

\newtheorem{thm}{Theorem}[section]
\newtheorem{lem}[thm]{Lemma}
\newtheorem{prop}[thm]{Proposition}

\newtheorem{rem}[thm]{Remark}

\newcommand{\mc}[1]{{\mathcal #1}}

\newcommand{\mb}[1]{{\mathbf #1}}
\newcommand{\bb}[1]{{\mathbb #1}}
\newcommand{\reff}[1]{(\ref{#1})}

\newlength{\defbaselineskip}
\newcommand{\setlinespacing}[1]%
           {\setlength{\baselineskip}{#1 \defbaselineskip}}

\title[Nonequilibrium Colour Fluctuations]{Nonequilibrium Density Fluctuations for the
Zero Range Process with Colour }
\author{Hanna K. Jankowski}

\address{\noindent Department of Statistics, University of
Washington
\newline
e-mail:  \rm \texttt{hanna@stat.washington.edu} }

\date{July 7, 2006}

\begin{document}

\maketitle

\begin{abstract}
We examine the fluctuations of the empirical density measure for
the colour version of the symmetric nearest neighbour zero range
particle systems in dimension one.  We show that the weak limit of
these fluctuations is the solution of a system of coupled
generalized Ornstein-Uhlenbeck processes.  We also discuss how
this result may be used to prove a central limit theorem for the
tagged particle on the level of finite dimensional distributions,
and identify the limiting variance. This is the central limit
theorem associated to propagation of chaos for this interacting
particle system.
\end{abstract}

\section{Introduction}

The zero range particle system describes a class of microscopic
models. It was first introduced by Spitzer in 1970~\cite{Sp} as an
example of an interacting particle system, and has since been
studied at length. One of its greatest advantages is its
mathematical tractability. It has found widespread application in
the modelling of nonequilibrium phenomena. Examples of these
include models of sandpile dynamics and other flow mechanisms, as
well as the repton model of gel electrophoresis~\cite{E}.

In the zero range system we consider particles which begin from a
random configuration and move around the discrete circle
$\bb{Z}/N\bb{Z}$ following these dynamics: each particle waits an
exponential amount of time to jump, it then chooses one of its
nearest neighbours with equal probability.  The exponential rate
for the first particle to leave a particular site is a function
$c(\cdot)$ of the number of particles at this site, and it is for
this reason that the system is called the zero range model.
Configurations of the system are denoted by $\eta$; if we are
currently at site $x \in \bb{Z}/N\bb{Z}$ then the number of
particles at the site is denoted by $\eta(x)$. Next, we
differentiate between particles by assigning to each one of $k$
colours.  The dynamics for each particle are the same as
previously, however, in the  colour process we keep track of the
number of particles of each colour at all sites $x$ of the
discrete circle, and we denote this as $\eta_i(x)$. If we ignore
the colour of each particle then we obtain the colour-blind model,
which simply keeps track of the total number of particles at each
site $x$.  We will say that the colour version \emph{contracts} to
the colour-blind process.

The dynamics we have described conserve the total number of
particles of each colour.  Furthermore, the configurations have a
family of invariant measures, indexed by the average density
vector $\bm \rho$.  These are usually referred to as the grand
canonical measures, and we denote them here by $\mu_{\bm\rho}$.
For the colour-blind model these contract to the measures
$\mu_\rho$ indexed by the total particle density.

The above construction creates the model on a microscopic level.
We, however, are interested in the system from a macroscopic
viewpoint.  In order to achieve this we re-scale space by
$\frac{1}{N}$.  We then need to speed up time by $N^2$ so as to
obtain a non-trivial system evolution.  This is the standard
diffusive space-time re-scaling.  We will denote the
configurations of the re-scaled system by $\eta^N$ and $\eta_i^N$
for the colour-blind and $k$-colour systems respectively.  After
re-scaling, the particles move around the discrete subset $\bb
T^N=\{0,\ldots,(N-1)/N\}$ of the unit circle $\bb T$.

We first wish to understand the evolution of the macroscopic
densities for the colour process.  To do this, we study the
asymptotic behaviour of the empirical densities
\begin{eqnarray*}
&\Pi_{i,t}^N = \frac{1}{N}\sum_{x\in \bb T^N} \eta_{i,t}^N(x).&
\end{eqnarray*}
 We
assume that the system is started in such a way so that the
initial empirical densities correspond to some fixed macroscopic
densities
\begin{eqnarray}\label{inicondintro}
&\lim_{N\rightarrow\infty}\Pi_{i,0}^N = \rho_{0}^i(x)dx,&
\end{eqnarray}
in the sense that this limit exists weakly on $\bb T$.

In the colour-blind model it is well known that (see for example
\cite{K-L}) that the macroscopic densities $\rho_t(x)$ evolve
according to the equation
\begin{eqnarray}\label{HSL1intro}
&\partial_t\rho = \frac{1}{2}\nabla D(\rho) \nabla \rho,& \
\rho_t|_{t=0}=\rho_0
\end{eqnarray}
with $\rho_0=\sum \rho^i_{0}$ given above.  We call $D(\rho)$ the
\emph{bulk} diffusion coefficient.  Results of this type are known
as hydrodynamic scaling limits in the literature.

It is also known that the empirical measures $\Pi_{i,t}^N$ satisfy
a law of large numbers and converge weakly to a non-random limit
$\rho_{t}^i(x)dx$ (see \cite{Ga-J-L} for a discussion).  The
vector $\bm \rho_t$, describing the limit, is the unique solution
of the coupled non-linear differential equation
\begin{eqnarray}\label{HSL2intro}
\partial_t \bm \rho = \frac{1}{2}\nabla \bm D_k(\bm \rho)\nabla \bm \rho,
\end{eqnarray}
with boundary condition $\bm\rho_0=\{\rho^i_{0}\}_{i=1}^k$. Here
$\bm D_k(\bm \rho)$ is the colour diffusion coefficient and is
given by the formula
\begin{eqnarray}\label{kbulkintro}
\bm D_k(\bm \rho)=S(\rho)I +
\frac{(D-S)(\rho)}{\rho}\left[\begin{array}{ccc}
          \rho^1 & \hdots & \rho^1 \\
          \vdots &        & \vdots \\
          \rho^k & \hdots & \rho^k \\
          \end{array}
        \right],
\end{eqnarray}
where $I$ is the $k\times k$ identity matrix and $S(\rho)$ is the
\emph{self} diffusion coefficient which we now describe.

Consider the zero range particle system where we have singled out,
or \emph{tagged}, one of the particles.  The remaining particles
are combined together to define a random environment of the tagged
particle.  This new system is a Markov process with identifiable
invariant measures. Using the methods of \cite{K-V} it is
straightforward to show that this tagged particle started in
equilibrium will converge to a diffusion with generator
$$\mc L = \frac{1}{2}S(\rho)\triangle.$$  $S(\rho)$ is
called the self diffusion coefficient.  One of the nice features
of the zero range model is that we can calculate $S(\rho)$
explicitly. Additionally, in the zero range case, we may express
the difference between the bulk and self diffusion coefficients,
\linebreak $D-S,$ as $S'(\rho)\rho$.  We calculate these
coefficients explicitly for the zero range model in
Remark~\ref{calcS}.  We also use the notation $\chi(\rho)$ to
denote the static compressibility.  For the zero range model this
is the same as the variance of the variable $\eta(x)$ under the
invariant measure $\mu_\rho$, and the following identity links the
static compressibility with the bulk and self diffusion
coefficients:
$$S(\rho)\,\rho=\chi(\rho)\cdot D (\rho).$$
It is important to note that this is a particular property which
holds for the zero range process, but not for general systems.

The hydrodynamic scaling limit for the colour densities of an
interacting particle system was first studied in \cite{Q} for the
symmetric simple exclusion process. This is a more difficult
problem as the simple exclusion colour process is of
\emph{non-gradient}~type.

 In this paper we study the central limit
theorem associated to the law of large numbers described above.
That is, we consider the fields defined by
\begin{eqnarray}\label{defYi}
<Y_{i,t}^N,f>=\sqrt{N}\left\{\frac{1}{N}\sum_{x\in\bb{T}_N}f(x)[\eta^N_{i,t}(x)-\rho_{t}^i(x)]\right\}.
\end{eqnarray}
The main result is the following theorem, which identifies the
asymptotic behaviour of $\bm Y^N_t.$ Let $\bm{\rho}_t$ be the
solution of \reff{HSL2intro}. Denote by $P_N$ the measure on
$D\{[0,T],H_{-4}^k\}$ induced by the stochastic field $\bm Y^N$.

\begin{thm}\label{thm:main} Under the assumptions given in
Section 3.3,  $P_N$ converges weakly in $D\{[0,T],H_{-4}^2\}$ to a
coupled generalized Ornstein-Uhlenbeck process characterised by
the coupled generalized stochastic differential equation
\begin{align}\label{resultintro} d\bm{Y}_t=\frac{1}{2}\{\triangle
\bm D_k(\bm \rho_t)\}\bm{Y}_t dt +\{ \nabla \bm A_k ^{1/2}(\bm
\rho_t )\}d\bm{W}_t
\end{align}
with initial condition $\bm{Y}_0 \sim P_0$. In the above $\bm D_k$
is the colour diffusion coefficient described previously, $\bm
A_k$ is given by
\begin{eqnarray}
\bm A_k(\bm \rho)= \frac{[\chi D](\rho)}{\rho}
\mbox{diag}\{\rho^1, \ldots, \rho_k\},\label{def:Ak}
\end{eqnarray}
where $\chi(\rho)$ is the static compressibility, and $D$ is again
the bulk diffusion coefficient. $\bm{W}_t$ is the $k$-dimensional
Gaussian random field defined as a vector composed of $k$
independent copies of the Gaussian random field $ W_t^i$, with
covariance
$$E[\langle W^i_{t},f \rangle \langle  W^i_{s},g \rangle] = \min(s,t) <f,g>,$$
for $i=1,\ldots,k.$
\end{thm}
The emergence of the colour diffusion coefficient in
\reff{resultintro} is not surprising, as the equation may be
obtained by formally linearizing the diffusion equation satisfied
by the macroscopic density profiles.  This notion is made precise
in a result called the Boltzmann-Gibbs principle.

Equation \reff{resultintro} contracts to the colour-blind density
fluctuation field which we may write in the form
\begin{eqnarray}\label{bulkresultintro}
dY_t = \frac{1}{2} \triangle D(\rho_t) Y_t dt + \nabla\sqrt{[\chi
D](\rho_t)} dW_t.
\end{eqnarray}
We next consider the deviations of the colour fields away from
$Y_t$.  To this end, define $U_t^i = Y_t^i -
\frac{\rho_t^i}{\rho_t} Y_t$. The generalized SDE satisfied by
$\bm U_t$ may now be written as
\begin{eqnarray}\label{devresultintro}
d\bm U_t = \frac{1}{2}\triangle S(\rho_t)\bm U_t dt + \nabla
\sqrt{S(\rho_t)} d\bm Z_t,
\end{eqnarray}
where $\bm Z_t$ is a Gaussian process \emph{independent} of $W_t$,
such that $1 \cdot \bm Z_t=0$.  In terms of the original
formulation we may express $dZ_t^i$ as
$\sqrt{\rho_t^i}dW_t^i-\frac{\rho_t^i}{\rho_t}\sum_{i=1}^k\sqrt{\rho_t^i}dW_t^i$,
which is independent of
$\sum_{i=1}^k\sqrt{\rho_t^i}dW_t^i=\sqrt{\rho_t}dW_t$ by direct
calculation of the covariances.  Thus we have that $\bm Z_t$ has
covariance structure given by
\begin{eqnarray}\label{Zcov}
E[\langle Z^i_{t},f \rangle \langle Z^j_{s},g \rangle] =
\int_0^{\min(s,t)}\int
f(x)g(x)\rho_u^i(x)\left\{\delta_{i,j}-\frac{\rho_u^j(x)}{\rho_u(x)}\right\}dxdu.
\end{eqnarray}
Given \reff{devresultintro} and \reff{bulkresultintro}, one may
recover \reff{resultintro}.

It is our belief that for general interacting particle systems,
such as symmetric simple exclusion, the colour density
fluctuations will be described by equations \reff{bulkresultintro}
and \reff{devresultintro}.   It is only the special relation
$[\chi D]\rho = S(\rho)\rho$ which allows the formula
\reff{resultintro} to emerge for the zero range process

Using equations \reff{resultintro} and \reff{bulkresultintro}, we
may also write the colour density fields as
$$dY^i_t = \frac{1}{2} \triangle S(\rho_t) Y^i_t dt+\frac{1}{2} \triangle \frac{[D-S](\rho_t)}{\rho_t} \rho_{t}^i Y_t dt + \nabla\sqrt{\frac{[\chi\cdot D](\rho_t)}{\rho_t}\rho^i_{t}} dW^i_t,$$
for $i=1\ldots k.$  The advantage of this formula is that we can
see explicitly how the density fluctuation field for any colour
interacts with the environment created by the colour-blind
process.  If our system was formed instead by non-interacting
random walks the middle term in the above formula would be equal
to zero.

In the literature, results on the fluctuations of the hydrodynamic
scaling limit are known as \emph{non-equilibrium} density
fluctuations. When the system is started in an equilibrium
measure, the terminology equilibrium density fluctuations is used.
In fact, the term ``density" is often omitted in the discussion;
however, we introduce it to emphasize the difference between these
central limit theorems and fluctuations for the tagged particles.

Next consider the empirical measure for the particle paths
themselves.  We denote the paths of the diffusively re-scaled
particles as $X_i^N(\cdot)$, each taking values in $D([0,T],\bb
T)$.  For each $N$ we have $n$ particles and we assume that
\begin{align}\label{littlen}
\lim_{N\rightarrow\infty}\frac{n}{N}=\bar{\rho}
\end{align}
exists.  The $n$ particles are initially randomly distributed on
the unit circle. Define the particle path empirical measure as
\begin{eqnarray}\label{pathemp}
&\Pi_n=\frac{1}{n}\sum_{i=1}^{n} \delta_{X^N_i(\cdot)}.&
\end{eqnarray}
 $\Pi_n$ is
a random variable taking values in $M_1\{D([0,T],\bb T)\}$, the
space of probability measures on $D([0,T],\bb T)$ endowed with the
weak topology.  That is, a realization of $\Pi_n$ is a measure
obtained by assigning mass $\frac{1}{n}$ to each observed particle
path.

A known result is that $\Pi_n$ satisfies a law of large numbers
\begin{eqnarray}\label{llnintro}
&\Pi_n \rightarrow \Pi,&
\end{eqnarray}
where the limit $\Pi$ is a non-random element in the space
$M_1\{D([0,T],\bb T)\}$. In fact, $\Pi$ is the measure of a
diffusion with generator
\begin{eqnarray}\label{Pgenintro}
\frac{1}{2}S(\rho_t)\triangle,
\end{eqnarray}
 where $\rho_t$ is the previously
discussed macroscopic density and $S$ is the self diffusion
coefficient.  The law of large numbers in \reff{llnintro} implies
that a randomly chosen tagged particle has the asymptotic
distribution $\Pi$. It also implies that the particle becomes
independent of the environment created by the system.     This
type of law of large numbers is often called \emph{propagation of
chaos}  (see, for example, \cite{Sz89}).Propagation of chaos
implies a central limit theorem for the position of a randomly
chosen tagged particle.

Propagation of chaos was proved for the symmetric simple exclusion
process in \cite{Q,R}.  First, the unique limit of $\Pi_n$ is
identified using the colour hydrodynamic scaling limit \cite{Q}.
The fluctuation results of this paper imply that the hydrodynamic
scaling limit for the zero range process satisfies equation
\reff{HSL2intro}. Tightness in the space $M_1\{D([0,T],\bb T)\}$
would complete the result.  This last fact was proved in \cite{R}
for symmetric exclusion, and the method presented there may be
applied to the zero range process. Large deviations associated to
the above law of large numbers were also studied for the symmetric
simple exclusion process in $d\geq 3$ \cite{Q-R-V}.

The next goal is to study the associated central limit theorems.
That is, we wish to understand the limiting behaviour of
$$\Gamma_n = \sqrt{n}(\Pi_n-\Pi).$$
In the case of independent random walks it is easy to see that the
limiting distribution of $\Gamma_n$ is a mean zero Gaussian random
field with the identity variance operator.  When there is
dependence between the particles the limit remains Gaussian, but
we expect an additional source of variation caused by the
interaction.  In Section \ref{sec:full} we discuss how one may
obtain the above limit for the zero range process as a direct
consequence of our non-equilibrium colour density fluctuation
result.  The type of convergence is convergence in finite
dimensional distributions. We also calculate explicitly the
variance for a class of test functions, and provide a formula for
the inverse of the variance in the general case.

Hydrodynamic scaling limits and associated fluctuations are both
subjects of much interest in the literature.  Scaling limits and
equilibrium fluctuations are well understood for many models (see
for example \cite{K-L, ML,Rost,Lu, GKL}). However, there is still
much work to be done to fully understand non-equilibrium
fluctuations.  To our knowledge there are currently only partial
results for gradient models, and no known results in the
non-gradient case. We adopt here Chang and Yau's \cite{Y92} proof
for the Ginzburg-Landau model in dimension one to our zero range
models.  The proof relies on knowing the logarithmic Sobolev
inequality for the inhomogeneous zero range process \cite{me}.  We
apply this inequality in the proof of $k$-colour non-equilibrium
density fluctuations. The main idea behind the argument is that we
may think of the zero range process for the $i$-th colour as a
process in a random environment, and our use of the logarithmic
Sobolev inequality for the inhomogeneous setting reflects this
notion.  In order to make use of this result we make an additional
technical assumption. For example, we could assume that the
particle jump rate, $c(k)$, is linear for large $k$ (see Section
\ref{assume} for a complete discussion of the assumptions). We do
not need this assumption elsewhere in this work, however, the full
result holds only in this setting. The remainder of the work holds
under the usual assumptions on the zero range process, namely
\begin{eqnarray}\label{LG}
&(LG) \hspace{1cm} \sup_k|c(k+1)-c(k)| <\infty&
\end{eqnarray}
as well as a weak monotonicity condition
\begin{eqnarray}\label{M}
&(M) \hspace{1cm} \inf_k \{c(k+k_0)-c(k)\} >0,&
\end{eqnarray}
for some integer $k_0$.

Propagation of chaos and the associated fluctuations have been
studied previously for systems with mean-field interactions,
~\cite{Sz84,Sz85,Sh-T,BA-B,T}. However, the interactions between
particles in the mean-field setting are weaker than for the
zero-range process. The method presented here to obtain the result
for finite dimensional distributions has not been previously
applied to show a fluctuation result.  It is based on the ideas
developed in \cite{Q, R, Q-R-V}.  As the colour version of the
symmetric simple exclusion process is of non-gradient type, at
this time we cannot prove tagged particle fluctuations for this
model.

The outline of this paper is as follows.    We begin with a
discussion of the relationship between nonequilibrium density
fluctuations for the colour version of the process and the central
limit theorem for the particle paths.  This is Section
\ref{sec:full}, where we also provide the formulae discussed
above. In Sections \ref{sec:main} and \ref{sec:prelim}, we define
notation, give a complete summary of the assumptions made, and
give some preliminary results which we will use in the proof of
the main result. Since the proof of Theorem \ref{thm:main} is
quite involved, we first explain the result for the colour-blind
model. This is done in Section \ref{sec:neqone}.  In Section
\ref{sec:neqtwo} we explain how to extend this to the colour
version of the model.

\section{Fluctuations of Particle Path Empirical
Measures\label{sec:full}}

A classical result from probability theory identifies the
fluctuations of the empirical measure of independent random
variables as a Gaussian field with identity covariance operator.
More precisely, let $X_1, X_2, \ldots$ be a collection of
independent and identically distributed random variables with
values in some Polish space $\mc X$ and common distribution $\pi$.
We define the empirical measure, $\hat\pi_n$, as the random
measure created by putting mass $\frac{1}{n}$ at each of the $n$
observed $X_i$ random variables:
\begin{eqnarray*}
&\hat\pi_n = \frac{1}{n}\sum_{i=1}^n \delta_{X_i}.&
\end{eqnarray*}
The strong law of large numbers implies that
$$\hat\pi_n \rightarrow \pi$$
in probability in the weak topology of probability measures on
$\mc X$. To study the fluctuations of this convergence we re-scale
the quantity $\hat\pi_n -\pi$ by $\sqrt{n}$. That is, we now
consider the quantity
$$\bb G_n=\sqrt{n}(\hat\pi_n - \pi).$$  A consequence of the
classical central limit theorem is that
$$\bb G_n \Rightarrow \bb G,$$
in the sense of finite dimensional distributions, where  $\bb G$
is the Gaussian random field with mean zero and identity
covariance operator.  That is, for each measurable function $f:\mc
X \rightarrow \bb R$ in $L^2(\pi)$ such that $\int f d\pi=0$ we
have
$$<\Gamma_n,f> \,\,\, \Rightarrow \,\,\, <\Gamma,f>,$$ where
$<\Gamma,f>$ is a real-valued Gaussian random variable with mean
zero and variance $\int f^2d\pi$.

The colour density fluctuations imply a similar result for the
fluctuations of the empirical measure where the random variables
$X_i$ are taken to be paths of the nearest neighbour zero-range
particles after diffusive re-scaling. In particular, we have a
system of $n$ interacting particles with trajectories taking
random values in $D([0,T], \bb T)$.  One would expect that these
fluctuations also converge to a Gaussian random field; however,
the random variables under study are no longer independent and
this nontrivial dependence structure would introduce an additional
correlation in the limiting covariance.

For each $N$ we have $n$ particles and we assume \reff{littlen}.
The $n$ particles are initially distributed on $\bb T$ and we
assume that this initial distribution satisfies both a law of
large numbers and a central limit theorem.  That is, we assume
that
\begin{eqnarray*}
&\lim_{n\rightarrow\infty}\frac{1}{n}\sum_{i=1}^n\delta_{X_i^N(0)}=\rho_0
(x)dx&
\end{eqnarray*}
exists weakly on $\bb T$.  Notice that this implies $\int
\rho_0(x)dx=1$.  We also assume that
\begin{eqnarray*}
&\lim_{n\rightarrow\infty} \sqrt{n}
\left(\frac{1}{n}\sum_{i=1}^n\delta_{X_i(0)}-\rho_0
(x)dx\right)=\Gamma_0&
\end{eqnarray*}
 exists, where $\Gamma_0$ is a random element
of $\mc H_{-4}$. These assumptions are implied by assumptions (D1)
and (F2) of Section \ref{assume}.

Define the empirical measure as in \reff{pathemp}. We next
describe a heuristic approach used to identify its limit $\Pi$.
Consider the time-marginals of the process $\Pi_n$.  Define
\begin{eqnarray*}
&\Pi_n(t) = \frac{1}{n} \sum_{i=1}^{n} \delta_{X^N_i(t)}.
\end{eqnarray*}
From Theorem \ref{HSL} we have that $\Pi_n(t)$ satisfies a law of
large numbers and converges weakly to the solution of
\reff{HSL1intro}. We would expect that a tagged particle started
out of equilibrium would also converge to a diffusion with the
self diffusion coefficient $S(\rho_t)$ and a drift term caused by
the evolution of the system towards equilibrium. That is, the
generator of this diffusion should be
\begin{align}\label{selfgentemp}
\mc L =\frac{1}{2} S(\rho_t)\triangle+ b\cdot \nabla.
\end{align}
We next identify the drift $b$ by noting that under the tagged
limit we would still expect the time marginals of $\Pi$ to evolve
according to the hydrodynamic diffusion equation \reff{HSL1intro}.
That is,
$$\mc L^\ast \rho_t = \frac{1}{2} \nabla \cdot D(\rho_t) \nabla \rho_t.$$
Equating with the forward equation from the generator, we obtain
\begin{align*}
b = \frac{1}{2 \rho_t}\{\nabla S(\rho_t)\rho_t-\nabla
D(\rho_t)\nabla\rho_t\}.
\end{align*}
In the zero range process we have $b\equiv0$, as can be shown
through direct computation of the above quantities. That is, for
the zero range process, the limiting distribution $\Pi$ is equal
to the measure $P$, where $P$ is the law concentrated on
$C([0,T],\bb T)$ of a drift-free diffusion process with generator
given by
\begin{align*}
\mc L =\frac{1}{2} S(\rho_t)\triangle.
\end{align*}

\bigskip

We wish to study the fluctuations of the previously described law
of large numbers.  That is, we center and re-scale as per usual,
and study the quantity
\begin{align}\label{fullfluctintro}
\Gamma_n = \sqrt{n}(\Pi_n-\Pi).
\end{align}
We have already identified $\Pi$ as the measure $P$ defined above.
Theorem \ref{thm:main} implies that $<\Gamma_n,F>$ converges to a
Gaussian random variable for smooth functions $F:C([0,T],\bb
T)\mapsto\bb R$ of the form
$$F(z(\cdot))=f(z(t_1),\ldots, z(t_k)),$$
where $t_1<\hdots<t_k$.  Functions of this form correspond to the
question ``\emph{what are the particles doing at
$\{t_1,\ldots,t_k\}$?}".  We answer the question by colouring the
particles according to their behaviour at these times, and hence,
we may use the fluctuation results  for the colour densities.

To describe this consider first the simplest case of
$$F(z(\cdot))=f(z_{t}).$$ Denote by $\mc P^D_{t,s}$ the semi-group associated to the
nonpositive operator $\frac{1}{2} D (\rho_t)\triangle$. The
results of Section \ref{sec:neqone} tell us that for this choice
of $F$, $<\Gamma_n, F>$ converges to the random variable described
by
\begin{eqnarray}\label{chfn:bulk}
E[e^{i<\Gamma, F>}]=\widehat{\pi}_0(\mc P^D_{t,0}f)
\exp\{-\frac{1}{2}E[\left\{\int_0^t\langle\mc
P^D_{t,s}f,dZ_s\rangle\right\}^2]\},
\end{eqnarray}
where $\widehat{\pi}_0(f)$  is the characteristic function of the
random variable $<\Gamma_0,f>.$  In the above $Z_t$ denotes the
generalized Gaussian process with covariance
\begin{align}\label{defZ}
E_P[\langle Z_t, f \rangle \langle Z_s, g
\rangle]=\int_0^{\min(s,t)}<\nabla f, [\chi D](\rho_u)\nabla g>du.
\end{align}
Indeed, this is simply a restatement of Theorem \ref{uniqueone}.

Consider next the case of
\begin{eqnarray}\label{chooseF2}
F(z(\cdot))=\bb I[z(t_1)\in B_1]f(t),
\end{eqnarray}
where $t>t_1$ and $B_1$ is an open subset of $\bb T$. Following
\cite{R}, we define
\begin{eqnarray}\label{reza}
\begin{array}{ccc}
\rho^1_{t_1}(x)&=&\rho_{t_1}(x)\bb I [x\in B_1]\\
\rho^2_{t_1}(x)&=&\rho_{t_1}(x)-\rho^1_{t_1}(x).
\end{array}
\end{eqnarray}
Let $\mc P^S_{t,s}$ denote the semi-group associated to the
equation $\partial_t  h = \frac{1}{2} S(\rho_t)\triangle  h.$  We
also define $\mc A^1_t$ as the operator acting on functions $h$ as
$\frac{(D-S)(\rho_t)}{\rho_t} \rho_{1,t} \triangle h.$ For the
zero-range process this becomes $S'(\rho)\rho^1\triangle$. Lastly,
as in the above, we define $Z_t^i$ to be two independent
generalized Gaussian processes with covariance
$$E_P[\langle Z_t^i, h \rangle\langle Z_s^i, g \rangle]=\int_0^{\min(s,t)}<\nabla f, [\chi D](\rho_u)\frac{\rho^i_{u}}{\rho_u}\nabla g>du.$$

The results of Section \ref{sec:neqtwo} imply that for the choice
of $F$ given as in \reff{chooseF2}, $<\Gamma_n,F>$ converges
weakly to the random variable described by
\begin{eqnarray}
E[e^{i<\Gamma, F>}]&=&\widehat{\pi}_0\left(\mc P^D_{t_1,0}\bb
I_{B_1}\mc P^S_{t,t_1}f +
\int_{t_1}^t \mc P^D_{s,0}\mc A^1_s\mc P^S_{t,s}f ds\right) \label{ini}\\
&&\times \exp\{-\frac{1}{2}E\left[\left\{\int_{t_1}^t\langle\mc
P^S_{t,s}f, dZ_s^1\rangle
\right.\right.
+ \int_0^{t_1}\langle \mc P^D_{t_1,u}\bb I_{B_1}\mc
P^S_{t,t_1}f, dZ_u^1+dZ_u^2\rangle\notag\\
&&\hspace{4cm} \left.\left.+ \int_{t_1}^t\int_0^{s}\langle\mc
P^D_{s,u}\mc A^1_s\mc P^S_{t,s}f, dZ_u^1+dZ_u^2\rangle
ds\right\}^2\right]\}.\notag
\end{eqnarray}
Notice that in the above formula the distribution of $Z_t^1+Z_t^2$
is the same as the distribution of the field $Z_t$ with covariance
given in \reff{defZ}.  We chose to leave the formula in this
format to emphasize the relationship between the terms.  This
formula is obtained by formally solving the generalized
Ornstein-Uhlenbeck stochastic differential equation.  That we may
do this is made rigorous by Theorem \ref{uniqueone} together with
Theorem~\ref{uniquetwo}.

We may now repeat the above to obtain convergence for general
functions $F$ for any number of $k$. That is,
\begin{eqnarray}\label{fnform}
F(z(\cdot))=\Pi_{i=1}^k\bb I[z(t_i)\in B_i]f(z(t)).
\end{eqnarray}
The $t_i$ are fixed and increasing times in $[0,T]$ with
$t_k<t<T$, $B_k$ are fixed measurable subsets of $\bb T$, and $f$
is a smooth function. Notice that each step doubles the number of
colours used, and we have $2^k$ colours for the general case.

Except for the initial measure $\pi_0$ of $\Gamma_0$, our limits
are Gaussian. If we assume that $\pi_0$ is also a Gaussian measure
on $\mc H_{-4}$ then we would have that $\Gamma$ is truly a
Gaussian field. We do not make this assumption; however, from this
point on we will say that $\Gamma$ is a Gaussian random field and
omit the discrepancy caused by the initial measure.

It follows from the above that we have proved weak convergence of
$$<\Gamma_n,F> \Rightarrow <\Gamma,F>$$
in the sense of finite dimensional distributions for a large class
of functions on \linebreak $C([0,T],\bb T)$. We know that $\Gamma$
is a mean zero Gaussian random field. We would like to identify
explicitly the covariance operator of $\Gamma$ for a general
function $G$.

For a field $G:C([0,T],\bb T )\rightarrow \mathbb{R}$ define
$$q(t,x)=E_P[G(X(\cdot))|X_t=x]\rho_t(x),$$
where $\rho_t$ is again the solution to \reff{HSL1intro}. Let
$\gamma^*$ denote the quantity
\begin{align*}
\gamma^*(t,x)=\frac{\nabla^{-1}\left\{\mc G
(t,\rho_t(x);q(t,x))\right\}}{[\chi\cdot D ](\rho_t(x))},
\end{align*}
where $\mc G$ is the functional
$$\mc G(t,\rho,q)= \partial_tq-\frac{1}{2}\triangle [D(\rho)q].$$
Lastly, define the operator $\mc A^*G = \int_0^T
\gamma^*(t,X_t)dX_t$.

\begin{thm}\label{myconjv2}
For any mean-zero function $G$ in $C_b\left[C([0,T],\bb T),\bb
R\right]$ into the real numbers, $<\Gamma,G>$ is a real-valued
Gaussian random variable described by
\begin{align}\label{finalformv2:chfn}
E[e^{i<\Gamma,G>}]=\hat\pi_0(F_0)\ e^{-\frac{1}{2} \mc Q_D(G)}\
e^{-\frac{1}{2}\mc Q_S(G)}.
\end{align}
 $\hat\pi_0$ is the
characteristic function of the initial field $\Gamma_0$ and $F_0 =
E_P[G(X(\cdot))|X_0]$ is the projection of the function $G$ onto
the initial field. $\mc Q_D(G)=E_P[G \,\,\Theta_D^{-1}G]$ may be
identified through the operators $\Theta_D$ on a test function $F$
corresponding to the operator on a field $G$ given by the
quadratic form
\begin{eqnarray} \label{finalformv2:bulk}
 E_P[(\mc A^*G)^2].
\end{eqnarray}
Similarly, for $\mc Q_S(G)=E_P[G \,\,\Theta_S^{-1}G]$ may be
identified as the inverse of the quadratic form
\begin{eqnarray} \label{finalformv2:dev}
 E_P[(G-\mc A^*G)^2].
\end{eqnarray}
\end{thm}

\begin{rem} In the case when we are working with a system of
independent random walks (that is, we choose $c(k)=k$), the above
formula becomes $E_P[G^2]$ as expected.
\end{rem}

\begin{proof}
It is enough to establish equivalence with the limiting behaviour
of $<\Gamma, F>$ for functions $F$ of the form in \reff{fnform}.

For any function $G$, the behaviour of $<\Gamma, G>$ may be split
into three components.  First, we have the behaviour of $<\Gamma,
E_P[G(\cdot)|X_0]>$, which is independent and given by the initial
field measure $\pi_0$. Next, we consider the evolution of the
marginals
$$<\Gamma, E_P[G(\cdot)|X_t]-E_P[G(\cdot)|X_0]>,$$ whose variance
is described by the operator $\mc A^* G = \mc A^*
[E_P[G(\cdot)|X_t]-E_P[G(\cdot)|X_0]]$.  This is the same as
\reff{bulkresultintro} and by formula \reff{finalformv2:chfn} is
independent of what remains.  The remaining portion is the
behaviour of $<\Gamma, G>$ less the time evolution of \linebreak
$<\Gamma, E_P[G(\cdot)|X_t]-E_P[G(\cdot)|X_0]>.$  By
\reff{finalformv2:chfn} and \reff{finalformv2:dev} this is
Gaussian with identity variance operator under the measure $P$.
For functions of the form given in \reff{fnform}, this is simply
\reff{devresultintro}, which is independent of
\reff{bulkresultintro}.  That is, consider the case of
$F(z(\cdot))=f(z_T)\bb I_B(z_s)$.  Using \reff{reza} and the
propagation of chaos results of \cite{R},  we have that
$E_P[f(X_T)\bb
I_{B}(X_s)|X_t=x]\rho_t(x)=E_P[f(X_T)|X_t=x]\rho^1_t(x)$, as in
the latter case the indicator function is simply equal to one. It
thus remains to argue that the appropriate density in this case is
given by
$\rho^i_t(x)\left\{1-\frac{\rho^i_t(x)}{\rho_t(x)}\right\}$ to
match \reff{Zcov}. This follows from noting that
$\rho_t(x)=\sum_{i=1}^k \rho_t^i(x)\bb I [\mbox{particle is of
type $i$}]$, for any number of colours. Thus, we may identify
$\Theta_D^{-1}$ via the variance given in \reff{chfn:bulk}, and
$\Theta_S^{-1}$ is simply the identity on $G-\mc A^*G$.
\end{proof}

Theorem \ref{myconjv2} tells us that the limiting fluctuations of
propagation of chaos for our model separate into the bulk density
fluctuations and the \emph{independent} fluctuations of the
remaining deviations.   Fluctuations of this type have been
studied previously for systems with mean-field interactions
in~\cite{Sz84,Sz85,Sh-T,BA-B,T}, where the limiting behaviour is
quite different.  For these models, the variance operator  may be
written in the form $(I-\mc A)^2$, with $\mc A$ described
explicitly in terms of Malliavin derivatives of $P$. The physical
source of this form for the variance operator is most easily seen
in the coupling approach developed in \cite{Sz84}.  We also note
that the convergence obtained in ~\cite{Sz84,Sz85,Sh-T,BA-B,T} for
the fluctuation results is weak convergence in finite dimensional
distributions only.


\section{Notation and Assumptions.\label{sec:main}}


\subsection{Notation.}

We denote the set $\{0,1/N,\ldots, (N-1)/N\}$ as $\bb T^N$, with
addition defined as addition modulo one, and the continuous unit
circle as $\bb T$.

Given a subset $\Lambda$ of $\bb T^N$ or of $\bb Z /N \bb Z$ we
write $AV_{x \in \Lambda}f(x)$ to denote the average of the
function $f$ inside the box $\Lambda$, that is,
\begin{eqnarray*}
AV_{x \in \Lambda}f(x)&=&\frac{1}{|\Lambda|}\sum_{x \in
\Lambda}f(x).
\end{eqnarray*}

For a metric space $\mc X$, $C^p(\mc X)$ stands for the space of
real-valued functions on $\mc X $ with $p$ continuous derivatives,
and $C_b(\mc X)$ the space of bounded and continuous functions on
$\mc X$.

For $z$ in $\bb Z$ define $h_z: \bb T \mapsto \bb R$ by
$e_z(x)=\sqrt{2}\cos(2\pi z u)$ for $z$ positive,
$e_z(x)=\sqrt{2}\sin(2\pi z u)$ for $z$ negative, and
$e_0(x)\equiv 1$.  The collection $\{e_z\}$ is an orthonormal
basis of $L^2(\bb T)$.  $<\cdot,\cdot>$ denotes the inner product
of $L^2(\bb T)$. The functions $e_z$ are also eigenfunctions of
the operator $I-\triangle$ , with eigenvalues $\gamma_z=1+4\pi
z^2$.  We use the notation $\mc H_m$ to denote the Hilbert space
formed by taking the completion of $C^{\infty}(\bb T)$, the space
of infinitely differentiable real functions on $\bb T$, under the
inner product
$$<f,g>_m=<f,(I-\triangle)^{m}g>.$$  For each positive integer
$k$ we denote by $\mc H_{-m}$ the dual of $\mc H_{m}$ relative to
the inner product $<\cdot,\cdot>$.

We define the k-fold product of these space to be
$$\mc H_{-m}^k=\mc H_{-m}\times\mc H_{-m}\times\hdots\times\mc
H_{-m}.$$ That is, for $\bm f \in \mc H_{-m}^k$ we have $\bm f =
\{f_1, \ldots, f_k\}$ and each $f_i$ is an element of $\mc
H_{-m}.$  The norm $||\bm f||_{-m,k}$ is defined as
$$||f_1||_{-m}+||f_2||_{-m}+\hdots+||f_k||_{-m}.$$

We also define the norm $||f||_{*}$ for a function defined on the
discrete unit circle and satisfying the condition $\sum_{x \in \bb
T^N}f(x)=0$ as
$$||f||_{*}=\frac{1}{N}\sum_{x\in \bb T^N}f(x)[-\triangle_N^{-1}f](x),$$
where $\triangle_N$ denotes the discrete Laplacian
$$[\triangle_Nh](x)=N^2\{h(x+1)-2h(x)+h(x-1)\}.$$  We denote by
$<\cdot,\cdot>_*$ the associated inner product. The traditional
notation for this norm is also $||\cdot||_{-1}$, which we do not
use to avoid confusion with the Sobolev spaces defined above. We
also define the product norm to be
$$||\bm f||_{*,k}=||f_1||_{*}+||f_2||_{*}+\ldots+||f_k||_{*},$$
for $\bm f = \{f_1, \ldots, f_k\}$.


\subsection{The Model\label{ch:def}}  The class of zero range
particles we describe is the symmetric nearest neighbour zero
range interacting particle system in dimension one.

\medskip

\noindent \textbf{The Colour-Blind Process.}  First consider the
evolution of the number of particles at each site. If a particle
moves from site $x$ to site $y$ the configuration $\eta$ changes
to $\eta^{x,y}$ where
$$(\eta^{x,y})(z) \; =\; \left\{
\begin{array}{ll}
\eta (z) - 1  & \hbox{if $z=x$}, \\
\eta (z) + 1  & \hbox{if $z=y$}, \\
\eta (z)      & \hbox{otherwise}.\;
\end{array} \right.
$$
It does this at rate $\frac{1}{2}c(\eta(x))$.  The system is a
Markov process and we can write down its generator $L$ as
\begin{equation}\label{def:gen}
(L f)(\eta) = \frac{1}{2}\sum_{x\sim y} c(\eta(x))
[f(\eta^{x,y})-f(\eta)]
\end{equation}
where $x\sim y$ denotes nearest neighbours of $\bb Z /N\bb Z$.


For $\varphi>0$ define the partition function $Z(\varphi) =
\sum_{k \ge 0} \frac{\varphi^k}{c(k)!}.$  Under the assumptions
(LG) and (M), the radius of convergence for $Z(\cdot)$ is
infinite. Fix $0<\varphi<\infty$ and denote by $\mu_\varphi$ the
product measure on $\mc X$ with marginals $$
\mu_\varphi(\eta(x)=k) = \frac{\varphi^k}{Z(\varphi)\,c(k)!}$$
where $c(k)! = \prod_{1 \le m \le k} c(m)$ if $k>0$ and $c(0)!=1$.
The family of measures $\mu_\varphi$ indexed by $\varphi$ are
stationary and reversible for the zero-range process. If we define
the product measures on the infinite lattice $\bb Z $, and not on
$\bb T$, these measures represent the full set of extremal
invariant measures for the system \cite{A}.

The measures $\mu_{\varphi}$ are referred to as the grand
canonical measures in the physics literature.  They are not
ergodic, however, as we have already pointed out that the system
evolution preserves the total number of particles.  If we
condition on the average number of particles, $y$, we do obtain
stationary, reversible, and ergodic measures:
$\mu_{\varphi}(\eta|\bar{\eta}=y)$.  Equivalently, we may also
condition on the total number of particles.  These conditional
measures are called the canonical ensembles.

To emphasize that these measures are defined on the discrete unit
circle of size $N$ we will use the notation $\nu_{N,y}$.  That is,
$\nu_{N,y}=\mu_{\varphi}(\eta|\bar{\eta}=y).$  If we are
considering the measures for configurations restricted to a subset
of the circle, and this subset is of size $K$, we will use the
notation $\nu_{K,y}$.  Because of the homogeneity of the system
there is no ambiguity in the notation.

Let $\rho(\varphi)=E_{\mu_{\varphi}}[\eta(x)]$ denote the average
density of particles.  By straightforward computation we obtain
the following identities
\begin{eqnarray}
\rho(\varphi)&=&\varphi\,Z^{\prime}(\varphi)\,Z^{-1}(\varphi),\\
\rho'(\varphi)&=&\varphi^{-1}\,E_{\mu_\varphi}[(\eta(x)-\rho(\varphi))^2]>0.
\end{eqnarray}
As the variance must be strictly positive for $\varphi>0$, we
conclude that $\rho(\varphi)$ is a strictly increasing function
and hence invertible.  Because $\rho$ has the natural
interpretation of the density, we will fix $\rho$ and think of
$\varphi$ as $\varphi(\rho)$. We shall then index the invariant
measure by $\rho$: $\mu_{\rho}$. Notice also that
$E_{\mu_{\rho}}[c(\eta(x))]=\varphi(\rho)$.

The static compressibility mentioned in the introduction is in
general defined as
\begin{eqnarray*}
&\sum_{x\in\bb Z/N\bb Z}E_{\mu_{\varphi}}[\eta(x);\eta(0)],&
\end{eqnarray*}
where we use the notation $E[f;g]$ to denote the covariance of the
functions $f$ and $g$.

We are now in the position where we may express the self and bulk
diffusion coefficients, as well as the static compressibility
explicitly.

\begin{rem}\label{calcS}
For the zero range process we have:
\begin{align*}
S(\rho)=\frac{\varphi(\rho)}{\rho}, \ \ D(\rho)=\varphi'(\rho), \
\ \mbox{and} \ \ \chi(\rho)=\frac{\varphi(\rho)}{\varphi'(\rho)}.
\end{align*}
\end{rem}

\noindent This follows from the above discussion of the invariant
measure $\mu_{\rho}$.
\medskip

We next define the Dirichlet form,
$D_{\mu_{\rho}}(f)=E_{\mu_\rho}[f(-L)f]$. Since $\mu_\rho$ is
reversible for the dynamics, this is equivalent to
\begin{align}\label{id:dirform}
D_{\mu_{\rho}}(f)=\frac{1}{2}\sum_{x\sim
y}E_{\mu_\rho}[c(\eta(x))[f(\eta^{x,y})-f(\eta)]^2].
\end{align}
The same identity holds when the expectation is taken with respect
to the canonical ensembles, in which case we denote the Dirichlet
form as $D_{\nu_{N,y}}(f)$. We may also restrict the Dirichlet
form to a subset of the discrete circle, $\Lambda_K$, where
$|\Lambda|=K$.  Here we again use the notation $D_{\nu_{K,y}}(f)$.
\begin{align}\label{id:dirform}
D_{\nu_{K,y}}(f)=\frac{1}{2}\sum_{x\sim y \in \Lambda_K
}E_{\nu_{K,y}}[c(\eta(x))[f(\eta^{x,y})-f(\eta)]^2].
\end{align}

\medskip

\noindent \textbf{The Colour Process.} To simplify notation we
define the model for the case when there are only two colours. The
extension to general $k$ is immediate.

Imagine that the zero range process is made up of two different
colours of particles. The two types of particles are mechanically
identical to the regular zero range process particles, but the
two-colour process keeps track of the two types of particles as
the system evolves. That is, we are now studying the evolution of
the number of particles of each type at site $x$ in $\bb Z/ N\bb
Z$.  Its elements will typically be denoted by the pair of
configurations \boldmath$\eta$\unboldmath$= (\eta_1,\eta_2) \in
{\mathbb N}^{\bb Z/ N\bb Z}\times{\bb N}^{\bb Z/ N\bb Z}$. The
dynamics for each particle are the same as in the zero range
process. Thus, the first particle of colour $i$ to jump from site
$x$ does so at rate
$$c_i(\bm{\eta}(x))=\frac{\eta_i(x)\,\,c(\eta_1(x))+\eta_2(x))}{(\eta_1(x)+\eta_2(x))}=\eta_i(x)\frac{c(\eta(x))}{\eta(x)},$$
where $\eta(x)=\eta_1(x)+\eta_2(x)$. We write down the
infinitesimal generator for the two colour process
\begin{eqnarray} (L f)(\bm \eta) & =
& \frac{1}{2}\sum_{i=1}^2 \sum_{x \sim y} c_i(\bm{\eta}(x)) [f(\bm
\eta^{x,y}_i))-f(\bm \eta)].\label{def:zr2}
\end{eqnarray}
Here $\bm \eta_i^{x,y}$ denotes the configuration obtained from
$\bm \eta$ by moving one particle of colour $i$ from site $x$ to
site $y$. Note that if the function $f$ is ``blind" to the
particle colour, i.e. $f(\bm{\eta}) = f(\eta_1 + \eta_2)$, then
$Lf$ is equivalent to the generator for the previously defined
zero-range process. Because of this contraction, we shall not use
a different notation for the generators of the two processes.

We now define the grand canonical measures and the canonical
ensembles for this process.  Fix $\bm \varphi =
\{\varphi_1,\varphi_2\}>0$ and denote by $\mu_{\bm \varphi}$ the
product measure on $\mc X^2$ with marginals
$$ \mu_{\bm\varphi}(\eta_1(x)=k, \eta_2(x)=m) =
\frac{\varphi_1^k \ \varphi_2^m}{c(k+m)!}\ {k+m \choose k} \
\frac{1}{Z(\varphi_1+\varphi_2)}.$$  $Z(\varphi)$ is the partition
function defined previously. We also have that
\begin{eqnarray} \mu_{\bm \varphi}( \eta_1(x)=k |
\eta_2(x)=m ) &=&
\frac{\varphi_1^k}{c_m(k)!Z_m(\varphi_1)},\label{def:condmeas}
\end{eqnarray}
where $c_m(k)=\frac{kc(k+m)}{k+m}$ and $Z_m(\varphi_1)$ is the
associated partition function.  The family of measures indexed by
$\bm \varphi$ are stationary and reversible for the two-colour
zero range process.

Let $ \rho^i(\bm\varphi) = E_{\mu_{\bm \varphi}}[\eta_i(x)]$
denote the density of particles of the $i^{\mbox{\tiny{th}}}$
colour. Notice that
\begin{eqnarray*}
\rho^i&=&\frac{\varphi_i\,Z'(\varphi)}{Z(\varphi)}, \mbox{ and}\\
\partial_{\varphi_j}{\rho^i}&=&\varphi_j^{-1}\,E_{\mu_{\bm
\varphi}}[\eta_i(x)\eta_j(x)-\rho^i\rho^j].
\end{eqnarray*}
In particular, this implies that the Jacobian of the
transformation $\bm \varphi \mapsto \bm\rho = \{\rho^1,\rho^2\}$
has determinant strictly positive for $\bm \varphi > 0$. As
before, we choose to index the invariant measure $\mu_{\bm
\varphi}$ by the pair $(\rho^1,\rho^2)$, (that is, work with
$\mu=\mu_{\rho^1,\rho^2}$), where we now consider
$\varphi_i=\varphi_i(\bm\rho)$.  The function $\varphi_i$ can be
recovered through $\varphi_i=E_{\mu_{\bm
\varphi}}[c_i(\bm{\eta}(x)]=\rho^i\frac{Z'(\varphi)}{Z(\varphi)}=\frac{\rho^i}{\rho}\varphi$.
Since
$c(\eta_1(x)+\eta_2(x))=c_1(\bm{\eta}(x))+c_2(\bm{\eta}(x))$, we
obtain that $\varphi=\varphi_1+\varphi_2 = E_{\mu_{\bm
\varphi}}[c(\eta_1(x)+\eta_2(x))]$ for the colour-blind model.

 As
before we define the canonical ensembles to be the measures
$\mu_{\bm \varphi}$ conditioned on the average density of both
colours of particles, which we denote as $\bm y =\{y_1,y_2\}$,
$$\nu_{N,\bm y}(\bm{\eta})=\mu_{\bm \varphi}(\bm{\eta}|\overline{\eta}_i=y_i,i=1,2).$$
If we consider this measure only for configurations restricted to
a subset of the circle, where the subset is of size $K$, we will
use the notation $\nu_{K,\bm y}$.

The Dirichlet form for the two-colour version, $D_{\mu_{\bm
\rho}}(f)=E_{\mu_{\bm \rho}}[f(-L)f]$, is then
\begin{align}\label{id:dirformtype}
&\frac{1}{2}\sum_{x\sim y}E_{\mu_{\bm
\rho}}[c_1(\bm{\eta}(x))[f(\bm
\eta_1^{x,y})-f(\bm{\eta})]^2]\notag+\frac{1}{2}\sum_{x\sim
y}E_{\mu_{\bm \rho}}
[c_2(\bm{\eta}(x))[f(\bm \eta_2^{x,y})-f(\bm{\eta})]^2]\notag\\
&=D^1_{\mu_{\bm \rho}}(f)+D^2_{\mu_{\bm \rho}}(f).
\end{align}

The same identity holds for the canonical ensembles, in which case
we have the Dirichlet form
$$D_{\nu_{N,\bm y}}(f)=E_{\nu_{N,\bm y}}[f(-L)f]=D^1_{\nu_{N,\bm y}}(f)+D^2_{\nu_{N,\bm y}}(f).$$ We will use the notation
$D_{\nu_{K,\bm y}}(f)$ as in the previous section, when we have
restricted the Dirichlet form to a subset of size $K$.

\medskip

We denote the expectation with respect to the measure induced by
the zero range process started in equilibrium as $E_{EQ}$ and as
$E_{NEQ}$ if the process is started out of equilibrium.  Because
of the previously noted contraction, there is no contradiction in
using this notation for both the colour and colour-blind models.

\subsection{Assumptions\label{assume}}

We separate the assumption into several subsections.

\medskip

\noindent\textbf{On the Rate Function. }

\smallskip

\noindent We make the assume that the rate function $c(\cdot)$
satisfies (LG) and (M) of \reff{LG} and \reff{M}. Assumption (LG)
is necessary to ensure that the zero range process is well defined
on the infinite lattice \cite{A}. Condition (M) rules out the
cases, such as the queueing system corresponding to $c(k)=\mathbb
I (k\geq 1)$, where the spectral gap is known to depend on the
density of particles.

A key ingredient in the proof of nonequilibrium density
fluctuations is the logarithmic Sobolev inequality.  In order to
make use of this tool in the inhomogeneous zero range process we
make an additional assumption.  We make the assumption \emph{only}
to be able to quote this result, and it is not required in the
remainder of the work.

\bigskip

\noindent \textbf{(E).}  Recall the conditional grand canonical
measure from \reff{def:condmeas}.  We assume that for every
$N_0\in \bb N$ there exists a finite positive constant $C=C(N_0)$
such that
\begin{eqnarray}
\frac{1}{C}\leq \sqrt{r}\, \mu_{N, \rho^1=\frac{r}{|\Lambda|} }
\left(\left.\sum_{x\in \Lambda}\eta_1(x)=r\right|\eta_2\right)
\leq C,
\end{eqnarray}
holds \emph{uniformly} over all $r\in\bb N$ and configurations
$\eta$ on a subset $\Lambda\subset\bb Z$ such that
$|\Lambda|=N_0$.

Because of the contraction principle, it is sufficient to make
this assumption for only two colours.  Also, by symmetry, it
follows that this assumption also holds for any colour,
conditioning on the configuration of the remaining particles. This
is a very technical condition, and it is satisfied if instead we
make, for example, one of the following assumptions.
\begin{enumerate}
\item[\textbf{(E1).}] There exists a large constant
$K_0$, and a positive constant $\theta$ such that for all $k\geq
K_0$ the rate function $c$ satisfies $c(k)=\theta k$.
\item[\textbf{(E2).}]  There exists a large constant $K_0$, and two
positive constants $\theta_1$ and $\theta_2$ such that for all
$k\geq K_0$ the rate function $c$ satisfies for all $x$
\[ c(k)=\left \{ \begin{array} {ll}
                            \theta_1 \, k & \mbox{if $k$ is odd,}\\
                            \theta_2 \, k & \mbox{if $k$ is even.}
                            \end{array}
                   \right. \]
\end{enumerate}
Naturally, many other variations on these exist.  For more details
on this assumption see \cite{me}.

\begin{rem}\label{rem:linerate}
The assumptions (LG) and (M) imply that there exist finite
positive constants $c_1$ and $c_2$ such that for all $k$ $c_1k\leq
c(k)\leq c_2k.$
\end{rem}


\bigskip \noindent\textbf{On Initial Density. }

\smallskip

\noindent Let $\bm{\zeta}_0^N$ denote a sequence of probability
measures on the k-fold product of $\mathbb{N}^{\mathbb{T_N}}$. For
each fixed $N$, $\bm\zeta_0^N$ stands for the initial measure of
the colour zero range process $\bm\eta^N_t$ with $k$ colours.

\bigskip
\noindent\textbf{(D1).} We assume that for each fixed $k$ the $\bm
\zeta_0^N$ are associated to some fixed density profile $\bm
\rho_0:\bb T^k\rightarrow[0,1]$ (where $\bm\rho_0 =
\{\rho_{0}^i\}_{i=1,\ldots,k}$ ) in the sense that
\begin{eqnarray*}
\lim_{N\rightarrow\infty}E_{\bm \zeta_0^N}\left[\left|
\frac{1}{N}\sum_{x\in
\mathbb{T_N}}f(x/N)\eta_i^N(x)-\int_{\mathbb{T}}
f(x)\rho_{0}^i(x)dx)\right|\right],
\end{eqnarray*}
for any continuous function $f\in C(\bb T)$ and for each $i$ from
1 to k .
\bigskip

\noindent \textbf{(D2).} We also make the assumption that
\begin{eqnarray*}
\limsup_{N\rightarrow\infty}E_{\bm
\zeta_0^N}\left[AV_{x}\{\eta^N(x)\}^2\right]\leq C.
\end{eqnarray*}


\bigskip

\noindent\textbf{On Initial Entropy. }

\smallskip

\noindent \textbf{(H1).} We assume that the relative entropy of
the intitial measure grows at most linearly with respect to an
invariant measure.  That is, for any $k$, assume that there exists
$\bm \varphi*>0$ and a constant $C=C($k$)$ such that
\begin{eqnarray*}
H(\bm \zeta_0^N|\mu_{\bm \varphi^*})\leq CN.
\end{eqnarray*}

\begin{rem}\label{rem:ent*}  Using the entropy inequality we can show that if the
above bound is satisfied for one vector $\bm \varphi^*$ then it is
also satisfied for each $\bm \varphi>0$.
\end{rem}

\bigskip
 \noindent \textbf{On Initial Fluctuations. }

\smallskip

\noindent Let $\mb Y_0^N = \{Y_{0,1}^N, \ldots, Y_{0,k}^N\}$,
where each $Y_{0,i}^N$ denotes the initial density fluctuation
field
$$Y_{0,i}^N(x)=\sqrt{N}\{\eta_i^N(x)-\rho_0^i(x)\}.$$
Also let $\bar{\mb Y}_0^N$ denote the centered version of $\mb
Y_{0}^N$.  That is
\begin{eqnarray}\label{def:Ybar}
\bar{Y}_{0,i}^N=\sqrt{N}\{\eta_i^N(x)-\rho_0^i(x)-AV_{x}(\eta_i^N(x)-\rho_0^i(x))\}.
\end{eqnarray}
\bigskip

\noindent \textbf{(F1).}  We assume the following
\begin{enumerate}
\item[(a).] $$\limsup_{N\rightarrow\infty}E_{\bm \zeta_0^N}[N^{-1/2}||\bar{\mb Y}_0^N||_{*,k}^2]=0$$
\item[(b).] $$\limsup_{N\rightarrow\infty}E_{\bm \zeta_0^N}[N^{1/2}|\overline\eta_i-AV_{x\in \bb T^N}\rho_0^i|^2]=0, \mbox{ for all $i=1,\ldots,k$.}$$
\end{enumerate}

\bigskip \noindent \textbf{(F2).}  Let $M(\mc H_{-m}^k)$ denote the
space of measures on $\mc H_{-m}^k$.  Assume that the initial
fluctuation field has a weak limit in $M(\mc H_{-m}^k)$ for $m=4$
. That is, if we denote by $P^N_0$ the measure induced on $M(\mc
H_{-m}^k)$ by the initial field $\bm Y_t^N$, then there exists a
unique measure $P_0$ such that $P^N_0\Rightarrow P_0$ in the weak
topology of $M(\mc H_{-m}^k)$.


\section{Preliminary Results\label{sec:prelim}}

We discuss in this section some results used in the remainder of
this work.

\subsection{Connection with Inhomogeneous Zero Range \label{relateZR}}

The symmetric inhomogeneous zero range process has dynamics which
are the same as of the homogeneous process, except that the jump
rate for a particle at site $x$ now also depends on the site $x$.
That is, denote the rates at site $x$ as $c_x(\cdot)$, and
$\xi(x)$ the number of particles currently at $x$.  Then the first
particle jumps from $x$ at rate $c_x(\xi(x))$, and it jumps to one
of its nearest neighbours with equal probability. We assume that
these rates satisfy the bounds (LG) and (M) uniformly in the
environment.
\begin{eqnarray*}
(LG^*) &\,&
\sup_{x,k}|c_x(k+1)-c_x(k)|\leq \tilde a_1 <\infty \\
(M^*) &\,& \inf_{x,k}|c_x(k+\tilde k_0)-c_x(k)| \geq \tilde a_2>0
\end{eqnarray*}
for some constants $\tilde a_1, \tilde a_2, \tilde k_0$.  The
generator of this process is
\begin{equation}
(\tilde{L} f)(\eta) = \frac{1}{2}\sum_{x\sim y\in\mathbb Z_N^d}
c_x(\xi(x)) [f(\xi^{x,y})-f(\xi)]
\end{equation}
The process has invariant grand canonical measures which are
product measures with marginals
$$\tilde{\mu}_{\Lambda,\varphi}(k)=
\frac{\varphi^k}{c_x(k)!\,\,Z_x(\varphi)}$$ where $Z_x(\varphi)$
is the appropriate partition function.  The process is reversible
with respect to these measures, as well as the canonical
ensembles, obtained by conditioning on the average density of
particles. We denote the canonical ensembles on a subset $\Lambda$
as
$\tilde\nu_{\Lambda,y}(\xi)=\tilde\mu_{\Lambda,\varphi}(\xi|\overline\xi
= y )$.

Notice that the mean density varies depending on site in the
inhomogeneous setting. However, we may still consider $\varphi$ as
an invertible function of the \emph{overall} density
$\frac{1}{|\Lambda|}\sum_{x \in
\Lambda}E_{\tilde\mu_{\Lambda,\varphi}}[\eta(s)]$ for every fixed
$\Lambda$. Similarly, we define
$\sigma_x^2(\varphi)=E_{\tilde\mu_{\Lambda,\varphi}}[\eta(x);\eta(x)]$,
and
$$\sigma^2_\Lambda(\varphi)=\frac{1}{|\Lambda|}\sum_{x\in\Lambda}\sigma^2_x(\varphi).$$ The
Dirichlet form, with expectation taken with respect to the
canonical ensembles on $\Lambda_K(x)=\{y: |x-y|\leq K \}$, shall
be denoted as
\begin{eqnarray*}
&\tilde{D}_{\tilde\nu_{\Lambda_K(x),y}}(f)=\frac{1}{2}E_{\tilde\nu_{\Lambda_K(x),y}}\left[\sum_{x\sim
y \in \Lambda_K(x)} \tilde
c(\xi(x))(f(\xi^{x,y})-f(\xi))^2\right].&
\end{eqnarray*}

We make the additional assumption that the inhomogeneous zero
range satisfies the following property.  For every $N_0\geq 2$,
there exists a constant $C=C(N_0)$ such that
\begin{eqnarray} (E^*) \hspace{.7cm} \frac{1}{C} \leq \sqrt{r}\,
\mu_{\Lambda, \varphi(\frac{r}{|\Lambda|}) } \left(\sum_{x\in
\Lambda}\eta(x)=r\right) \leq C
\end{eqnarray}
uniformly in $r\in \bb N$, and sets $\Lambda$ for any
$N_0=|\Lambda|$.

\medskip

\noindent\textbf{Relationship with Colour Zero Range.}  Consider
the colour zero range process where all except the particles of
the first colour have been ``frozen".  By the contraction
principle, we may consider only the case of two colours without
loss of generality.

To this end, fix a configuration $\eta_2$ and let $\eta_1$ evolve
as though it was a single-colour zero range process which uses the
\emph{non-homogenous} rate function
\begin{eqnarray*}
c_x(\eta_1(x))=\frac{\eta_1(x)\,\,c(\eta_1(x)+\eta_2(x))}{\eta_1(x)+\eta_2(x)}.
\end{eqnarray*}
This is a particular example of a non-homogeneous zero range
process. The invariant measure for this process is clearly the
product measure with marginals given in \reff{def:condmeas},
 the two-colour process invariant measure conditioned on the
second colour configuration.  It is not difficult to show that
conditions (LG), (M) and (E) imply that the measures
$\tilde{\mu}_x(k)$ satisfy conditions $(LG^*), (M^*) \mbox{ and}
(E^*)$.

The Dirichlet form (with respect to the canonical ensemble) in
this case is simply 
\begin{eqnarray}\label{conddirform}
&D^{x,1}_{\nu_{K,\bm y}(\cdot|\
\eta_2)}(f)=E_{\nu_{\Lambda_K(x),\bm y
}(\cdot|\eta_2)}\left[\sum_{x\sim y \in
\Lambda_K(x)}c_x(\eta_1(x))(f(\bm \eta_1^{x,y})-f(\bm
\eta))^2\right].&
\end{eqnarray}

We shall use the results of \cite{me} on the inhomogeneous zero
range process to state some useful facts about the conditional
colour zero range.  In what follows we shall state the results for
colour 1 conditioning on colour 2, or equivalently, on all of the
other colours.  By symmetry, the results are valid for any colour.

\subsection{Moment Bounds.}

The following two lemmas are proved in ~\cite{L-Se-V}.

\begin{lem}
There exist constants $C_1$ and $C_2$, which depend only on the
values $a_1, a_2$ and $k_0$, such that
$$0<C_1<\frac{\sigma^2(\rho)}{\rho}<C_2<\infty.$$
\end{lem}

\begin{lem}\label{lem:momentbounds}
For all $k\geq1$, there exists a finite constant $C(k)$ such that
$$m_{2k}(\rho)\leq C(k)\,\,\sigma^{2k}(\rho)$$
for all $\rho\geq\rho_0>0,$ where $m_j(\rho)$ denotes the
$j^{\tiny{th}}$ moment of $\eta(x)$ under the distribution
$\mu_\rho$.
\end{lem}

\noindent Next we prove a simple but useful property of the
functions $\varphi(\rho)$ and
$S(\rho)=\frac{\varphi(\rho)}{\rho}$:

\begin{prop}[Lipschitz properties]\label{ch1:prop:lip}
The following functions are Lipschitz:
\begin{enumerate}
\item $\varphi(\rho)$
\item $\varphi_i(\bm{\rho})$ as a function of $\rho^i$
\item $S(\rho)$
\end{enumerate}
The first two functions are also strictly increasing.
\end{prop}

\begin{proof}
The proofs follow from direct calculations of the derivatives, as
well as previously computed bounds:
\begin{enumerate}
\item $\varphi'(\rho)=\frac{\varphi(\rho)}{\sigma^2(\rho)}$.
\item
$\partial_{\rho^1}\varphi_1(\bm{\rho})=\frac{\rho^2}{\rho}S(\rho)+\frac{\rho^1}{\rho}\varphi'(\rho)$.
\item First of all notice that if $c(k)=k$ then $S(\rho)\equiv1$ and there is nothing to do.  Otherwise,
$ S'(\rho)\rho = \varphi'(\rho)-S(\rho)$ which implies that
$S'(\rho)$ is bounded below as long as $\rho$ is bounded below. To
finish we need to examine the limit of $S'(\rho)$ as
$\rho\rightarrow 0$. This follows from noting that
\begin{align*}
\lim_{\rho\rightarrow 0} S'(\rho) &= -\lim_{\rho\rightarrow 0}
S(\rho)^2 \frac{\varphi}{\sigma^2}\times \lim_{\rho\rightarrow
0}\left\{\frac{Z''(\varphi)Z(\varphi)-Z'(\varphi)Z'(\varphi)}{Z^2(\varphi)}\right\}.
\end{align*}
\end{enumerate}
\end{proof}

\subsection{Spectral Inequalities.}

When the underlying measure $\mu$ is one of the canonical
ensembles $\nu_{K,y}$ we will write the entropy as $H_{K,y}(f)$
for ease of notation.  The following inequality was proved by Dai
Pra and Posta~\cite{DP-P}, \cite{DPPG} under assumption (LG) and
(M).

\begin{thm}[Logarithmic Sobolev Inequality]\label{lem:logsob}
There exists a constant $C_{LS}=C_{LS}(a_1,a_2,k_0)$ such that
$$H_{K,y}(f)\leq C_{LS}K^2 D_{\nu_{K,y}}(\sqrt{f})$$
holds for any $K,y$ and positive function $f$.
\end{thm}

If we make the additional assumption (E), in light of the
discussion of section \ref{relateZR} and the results of \cite{me},
we obtain a logarithmic Sobolev bound for the conditioned zero
range process as well.

\begin{thm}[Logarithmic Sobolev Inequality for Conditioned Zero Range]\label{thm:logsobnh}
Assume that conditions (LG), (M), and (E) hold. Then there exists
a constant \, $\tilde{C}_{LS}= \tilde{C}_{LS} (a_1,a_2, k_0)$ such
that
$$H(f|\nu_{K,\bm y}(\cdot|\eta_2))\leq \tilde{C}_{LS}\,K^2 \,D^{x,1}_{\nu_{K,\bm y}(\cdot|\eta_2)}(\sqrt{f})$$
holds for any $\Lambda_K(x)$, $y_1$, density $f$ on $\bb N
^{\Lambda_K(x)}$, and configuration of particles of the second
colour $\eta_2$.
\end{thm}

\subsection{Hydrodynamic Scaling Limits.}

\begin{thm}\label{HSL}
Under assumptions (LG), (M), (D1) and (D2) we have that for every
$t \leq T$, for every continuous function $g:\mathbb T^d
\rightarrow \mathbb R$ and for every $\delta > 0$,
$$\lim_{N \rightarrow \infty }
P_{\zeta^N_0}\left[\,\,\,\left|\frac{1}{N}\sum_{x \in \mathbb
Z/N\bb Z} g(x/N) \eta_t(x) - \int _{\mathbb T} g(\theta)
\rho_t(\theta)d\theta\right|
> \delta \right] = 0$$
where $\rho_t(\theta)$ is the unique weak solution of the
non-linear heat equation \reff{HSL1intro}. Here we use
$P_{\zeta^N_0}$ to denote the measure of the zero range process
with initial measure $\zeta_0^N$.
\end{thm}

\begin{rem}
It is possible to prove this under reduced assumptions on the rate
function.  It is sufficient that $c(k)\geq \theta k$ for all
$k\geq 0$ and some positive constant $\theta$.
\end{rem}
\noindent The proof of the above fact appears in Kipnis and Landim
(\cite{K-L}) and is based on the entropy method first developed by
Guo, Papanicolaou, and Varadhan \cite{Gu-P-V}.

\begin{rem}
Theorem \ref{thm:main} implies the hydrodynamic scaling limit
result for the colour version of the zero range process.  However,
it is possible to prove this result independently and under fewer
assumptions.  Most notably, assumption (E) is not necessary.  We
may use, for example, the nonhomogeneous spectral gap from
\cite{me}.
\end{rem}

\subsection{Uniform Local Limit Theorems.}

For $m\geq0$, denote by $H_m(x)$ the Hermite polynomial of degree
$m$:
\begin{align*}\label{subsec:hermite}
H_m(x)=(-1)^m \exp {(-\frac{x^2}{2})} \frac{d^m}{dx^m} \exp
{(-\frac{x^2}{2})}.
\end{align*}
Let $g_0(x)$ denote the density of a standard normal random
variable, and define for $j\geq1$
\begin{align}
g_j(x)=g_0(x) \sum_{} H_{j+2a}(x) \prod_{m=1}^j
\frac{1}{k_m!}\left(\frac{\kappa_{m+2}}{(m+2)!\sigma^{m+2}}\right)^{k_m}
\end{align}
where the sum is taken over all nonnegative integer solutions
$\{k_l\}_{l=1}^j$ and $a$ such that \linebreak
$k_1+2k_2+\hdots+jk_j=j$ and $k_1+k_2+\hdots+k_j=a$.

\begin{thm}\label{thm:localCLT}
For all $\rho^*>0$ and $J \in \mathbb{N}$, there exist finite
constants $E_0=E_0(\rho^*,J)$ and $A=A(\rho^*,J)$ such that
\begin{align*}
\left|\sqrt{N\sigma^2(\rho)}\mu_{\rho}\left[\sum_{x\in\Lambda}\eta(x)=N\rho+\sigma\sqrt{N}z\right]-\sum_{j=0}^{J-2}\frac{1}{N^{j/2}}g_j(z)\right|
\leq \frac{E_o}{(\sigma^2(\rho)N)^{(J-1)/2}}
\end{align*}
uniformly over $z$ and over all parameters $\rho^* \geq \rho \geq
A/N.$
\end{thm}

The above is an Edgeworth expansion for a lattice distribution,
valid uniformly for the family of measures $\mu_{\rho}$, the first
as  $A/N \leq \rho \leq \rho_0$. Its proof appears in
\cite{L-Se-V} We also make use of the same result proved in
\cite{me} for the inhomogeneous process.  We state the result for
the conditioned zero range process.  In the theorem we have that
$|\Lambda|=N$.

\begin{thm}\label{normalLLT}
For all $\rho^*>0$ and $J \in \mathbb{N}$, there exist finite
constants $E_0=E_0(\rho^*,J)$ and $A=A(\rho^*,J)$ such that
\begin{align*}
\left|\sqrt{N\sigma_\Lambda^2}\mu_{N,\bm\rho}\left[\left.\sum_{x\in\Lambda}\eta_1(x)=N\rho^1+\sigma\sqrt{N}z\right|\eta_2\right]-\sum_{j=0}^{J-2}\frac{1}{N^{j/2}}g_j(z)\right|
\leq \frac{E_0}{(\sigma^2_\Lambda N)^{(J-1)/2}}
\end{align*}
uniformly over $z$, configurations of $\eta_2$, and over all
parameters $A/N \leq \rho^1 \leq \rho^*.$
\end{thm}


\section{Density Fluctuations in Non-equilibrium\label{sec:neqone}}


Let $\rho_t$ be the unique solution to the partial differential
equation \reff{HSL1intro}.   Define the density fluctuation field
for the colour-blind model as:
\begin{align}\label{deffieldone}
<Y_t^N,f>\,=\sqrt{N}\{\frac{1}{N}\sum_{x\in\bb{T}_N}f(x)[\eta^N_t(x)-\rho_t(x)]\}.
\end{align}
Let $P_N$ denote the measure of on $D\{[0,T],H_{-4}\}$ induced by
the stochastic field $Y^N$ with initial measure $\zeta_0^N$.
\begin{thm}\label{thm:neqone}
Under the assumptions of Section \ref{assume}, except for (E),
 $P_N$ converges weakly in $D\{[0,T],H_{-4}\}$ to a
generalized Ornstein-Uhlenbeck process characterised by the
generalized stochastic differential equation
\begin{align}\label{GenSDE1}
dY_t = \frac{1}{2} \{\triangle \varphi'(\rho_t)\} Y_t +
\{\nabla\sqrt{\varphi(\rho_t)}\} dW_t,
\end{align}
where $W_t$ is the Gaussian random field $\langle W_t,f\rangle$
with covariance
$$COV[\langle W_{t},f \rangle \langle W_{s},g \rangle] = \min(s,t)
<f,g>.$$
\end{thm}
The process $W_t$ is also known as a  generalized Brownian motion
which is ``Brownian" in time and ``white" in space.

Before proceeding with the proof of the above we make the quick
remark that non-equilibrium density fluctuations have been proved
for a specific class of zero range processes in \cite{FPV}. The
class of processes in their work is very special; the rate
function they consider is $c(k)=\bb I[k\geq 1]$. In this very
unique case there exists a special equivalence between this zero
range process and the symmetric simple exclusion process.  In this
setting, the martingale equations are closed in the field $Y_t^N$
and no Boltzmann-Gibbs argument in necessary.

\subsection{Outline of Proof\label{outline}}

The proof is divided into several main steps.  We first identify
the drift and quadratic variation of the limiting field $Y_t$.  We
then show that there is only one measure which solves this
martingale problem, and that the sequence of measures $P^N$ is
tight.  We first turn our attention to the drift and quadratic
variation martingales. Suppose that $P_{N_k}$ is a convergent
subsequence of $P_N$, and for ease of notation, we denote it again
by $P_N$. Let $P$ denote its weak limit.

Identifying the drift and quadratic variation of $Y_t$ which
matches the equation \reff{GenSDE1} is equivalent to the statement
that for any test function $f$
\begin{eqnarray}\label{limitmart1}
&<M_t,f>\ =\ <Y_t,f>-<Y_0,f>-\int_0^t
<Y_s,\frac{1}{2}\varphi'(\rho_s(\cdot))\triangle f>ds,&
\end{eqnarray}
is a martingale under  $P$ with quadratic variation
given by
\begin{eqnarray}\label{limitqv1}
&E_P[(<M_t,f>-<M_s,f>)^2|\,\mathcal{F}_s]=\int_0^t \int
\varphi(\rho_u(x))[\nabla f(x)]^2 dxdu.&
\end{eqnarray}
We next
compare these to the martingales under the measures $P_N$.

By direct calculation we know that $M_t^N$ defined below is a
martingale under $P_N$ for any $N$.
\begin{eqnarray}
&<M_t^N,f>
      =<Y_t^N,f>-<Y_0^N,f>-\int_0^t<\tilde{Y}_s^N,\triangle_Nf>ds&\label{line:martcalc}
\end{eqnarray}
where $\tilde{Y}$ denotes the fluctuation field
$\sqrt{N}(c(\eta^N(x))-\varphi(\rho(x)))$.

\medskip

\noindent\textbf{Comment on notation:}  There is a slight abuse of
notation in the above as we use the same inner product notation
$<f,g>$ for functions $f,g$ defined on $\bb T$, where it is equal
to $\int f g dx$, and for functions $f,g$ defined on $\bb T^N$
where it is equal to $\frac{1}{N}\sum_{x\in \bb T^N}f(x)g(x)$.  We
do this only to simplify the notation.

\medskip

A similar calculation gives the quadratic variation of
$<M_t^N,f>$: under $P_N$, $N_t^N(f)$ defined below is a
martingale.
\begin{align}\label{line:qvcalc}
         &<M_t^N,f>^2-\frac{1}{2}\int_0^t AV_{x\in
       \bb{T}^N}\left[c(\eta^N_s(x))+c(\eta^N_s(x+\frac{1}{N}))\right][\nabla_Nf(x)]^2ds.&
\end{align}
It hence remains to show that the limits of $M_t^N$ and $N_t^N$
are consistent with \reff{limitmart1} and \reff{limitqv1}. This is
much easier for $N_t^N(f)$. The difficulty which arises in $M_t^N$
is that expression \reff{line:martcalc} is not closed in the field
$Y^N$. This will prove to be the main obstacle to overcome in the
proof. We start with the fact that $<M_t^N,f>$ is a martingale,
and hence we know that for all $N$
\begin{align}\label{id:martN}
&E_{P_N}[<M_t^N,f>U]=E_{P_N}[<M_s^N,f>U],&
\end{align}
for all bounded, $\mathcal{F}_s$-measurable $U$. Here, $\mc F_t$
denotes the $\sigma$-algebra on $D\{[0,T],H_{-4}\}$ generated by
$F_s(f)$ for  $s\leq t$ and  $f\in C^{\infty}(\bb T),$ where $F
\in D\{[0,T],H_{-4}\}.$ We need to show
\begin{align}\label{id:mart}
&E_{P}[<M_t,f>U]=E_{P}[<M_s,f>U],&
\end{align}
Comparing \reff{limitmart1} with \reff{line:martcalc} we see that
we need to replace the field $\tilde Y^N_t$ with
$\varphi'(\rho_t)Y_t^N$. This type of result is known as the
Boltzmann-Gibbs principle in the literature.

The remainder of this section will be divided in the following
manner. We begin by showing that every weak limit of the measures
$P_N$ solves the martingale problem.  This is subdivided into
identifying the asymptotic drift and then its quadratic variation.
We first consider the drift, and begin with the proof of the
Boltzmann-Gibbs principle.  We handle the quadratic variation in
Section \ref{section:QV}. In Section \ref{sec:uniqueone} we
discuss the theory of Holley and Stroock with states that there is
only one solution to the martingale problem.  Section
\ref{sec:tightone} is dedicated to tightness of the measures
$P^N$. These results in combination as outlined above prove
Theorem \ref{thm:neqone}.

\subsection{Identifying the Drift Martingale.}

In light of the preceding discussion in the introduction, to show
that under $P$ the drift martingale is identified through
\reff{limitmart1}, it remains to prove the following.

\begin{thm}[Boltzmann-Gibbs Principle.]\label{BG}
For functions $f \in C^1[\bb{T}]$
\begin{eqnarray}\label{BGeqn}
\lim_{N\rightarrow\infty}E_{NEQ}\left[\left|\int_0^t
<f,\sqrt{N}(c(\eta^N_s)-\varphi(\rho_s)-\varphi'(\rho_s)(\eta^N_s-\rho_s))>ds\right|\right]\rightarrow0.
\end{eqnarray}
\end{thm}
\noindent This section is dedicated to the proof of this result.

The first step of the proof is to replace $\rho_t(x)$, the
solution of \reff{HSL1intro} with the solution to the discretized
version
\begin{align} \label{discreteHSLone}
\partial_t\rho=\frac{1}{2}\triangle_N \varphi(\rho),
\end{align}
with initial conditions $\rho_t(x)|_{t=0}=\rho_0(x)$ for $x$ in
$\bb T_N$.  The difference between the two solutions is of order
$\frac{1}{N}$  and hence does not affect \reff{BGeqn}, \cite{RM}.
Because of this fact and in order to simplify notation, we
continue to denote the solution of \reff{discreteHSLone} as $\rho$
in the remainder of this section.

We begin by re-writing the field
$\sqrt{N}(c(\eta^N_t)-\varphi(\rho_t)-\varphi'(\rho_t)(\eta^N_t-\rho_t))$
as the sum of five separate parts:
\begin{eqnarray}
\Phi^N_1(t)&=&\sqrt{N}[\eta^N_t(x)-\bar{c}_K(\eta^N_t(x))-\varphi'(\rho_t(x))\{\eta^N_t(x)-m^K_t(x)\}]\notag\\
\Phi^N_2(t)&=&\sqrt{N}[\bar{c}_K(\eta^N_t(x))-\varphi(m^K_t(x))]\bb{I}[m^K_t(x)\leq R]\notag\\
\Phi^N_3(t)&=&\sqrt{N}[\varphi(m^K_t(x))-\varphi(\rho_t(x))-\varphi'(\rho_t(x))\{m^K_t(x)-\rho_t(x)\}]\bb{I}[m^K_t(x)\leq R]\notag\\
\Phi^N_4(t)&=&\sqrt{N}[-\varphi(\rho_t(x))-\varphi'(\rho_t(x))\{m^K_t(x)-\rho_t(x)\}]\bb{I}[m^K_t(x)> R]\notag\\
\Phi^N_5(t)&=&\sqrt{N}\bar{c}_K(x)\bb{I}[m^K_t(x)>
R].\label{BGdiv}
\end{eqnarray}
Here, $m^K_t(x)=AV_{|y-x|\leq K/N}\,\eta^N_t(y)$ and similarly,
$\bar{c}_K(x)=AV_{|y-x|\leq K/N}\,c(\eta^N_t(y))$.

To prove Theorem \reff{BG}, we thus need to show that
\begin{eqnarray}\label{BGtodo}
&E_{NEQ}[|\int_0^t<f,\Phi_i^N(s)>~ds|]\rightarrow0&
\end{eqnarray}
 for for each
$i$. We begin by proving several lemmas, which summarize the two
main components of the Boltzmann-Gibbs principle.  The argument is
essentially a Taylor argument, and the local equilibrium principle
gives us the first term in the expansion, while ``equivalence of
solutions" allows us to control the resulting error term.


\subsubsection{Local Equilibrium Principle.} The first lemma we need
proves that over boxes of microscopic size $K$ we are able to
replace the average of a function with the expectation of said
function at the average density $m^K(x)$. We take $K<<N$.  It
turns out in the end that the correct scaling for this result is
to take $K$ to be slightly smaller than $\sqrt{N}$. To handle this
we will write $K=l\sqrt{N}$, and we will let $l\rightarrow0$.

Our main purpose in proving this result is to handle $i=2$ in
\reff{BGtodo}, however, we state and prove the result with
slightly more generality. For a function $g:\bb{R}\mapsto\bb{R}$
define $V(x)=V_{(g,K,R)}(x)$ to be the function
\begin{align*}
V(x)=(\bar{g}_K(\eta^N(x))-\hat{g}_K(m^K(x)))\bb{I}[m^K(x)\leq R],
\end{align*}
where $\bar{g}_K(\eta^N(x))=AV_{|x-y|\leq K/N}g^N(\eta(y)),$ and
$\hat{g}_K(\rho)= E_{\mu_{\rho=m^K(x)}}[g(\eta^N(y))]].$
 If the configuration $\eta$ depends on time, then $V$ will
also depend on time, and we shall write this as $V_s$. Notice that
in the $E_{\mu_{\rho}}[g(\eta^N(y))]]$ we are considering
configurations $\eta$ on $\bb T^N$, and not on the discrete
circle.  However, this does not change the invariant measures, and
for this reason we do not alter the notation.

\begin{lem}[Local Equilibrium Principle] \label{lem:LEP}

Suppose that $g$ is a function satisfying $|g(x)|\leq C(1+x)$,
then
\begin{align*}
\lim_{l\rightarrow0}\lim_{N\rightarrow\infty}\sup_{||J||_{\infty}<1}N^{1/2}E_{NEQ}[|\int_0^t<J,V_s>ds|]\leq0.
\end{align*}
\end{lem}
\begin{proof}
Let $X_s$ denote $<J,V_s>$. By the entropy inequality and
assumption (H1), we have that the following bound for any positive
$\beta$,
\begin{align}
E_{NEQ}[N^{1/2}|\int_0^tX_sds|] \leq  \frac{2}{\beta N} \log
E_{EQ}\left[\exp\left\{\pm \beta N \int_0^t
N^{1/2}X_sds\right\}\right]+\frac{C}{\beta}.\label{line:lep1}
\end{align}
It follows from the Feynman-Kac formula that
\begin{align}
E_{EQ}[\exp\{\pm \beta N \int_0^t N^{1/2}X_sds\}] \leq
\exp\{\int_0^t \Gamma^N_sds\}\label{line:lep2},
\end{align}
where $\Gamma^N_s$ is the largest eigenvalue of $\pm \beta
N^{3/2}X_s + L_N$. We next reduce the span of the Dirichlet form
$D_{\mu_\rho}(\sqrt f)$. We have that
\begin{align*}
D_{\mu_\rho}(\sqrt f)               & \geq \frac{N}{2K+1} \inf_{x}
\sum_{|x-z|\leq \frac{K}{N}}
\sum_{|z-y|=\frac{1}{N}}E_{\mu_{\rho}}[c(\eta^N(x))\{\nabla_{x,y}\sqrt{f}\}^2]
\end{align*}
where $\nabla_{x,y}f=f(\eta^{x,y})-f(\eta)$. We denote the right
hand side of the last line above by $\frac{N}{2K+1} \inf_{x}
D_{\mu_{\rho}}^{K,x} $. We may hence bound $\Gamma^N_s$ by
\begin{align}
N^{3/2} \sup_{f}
\sup_{|\alpha|\leq\beta||J||_{\infty}}\{E_{\mu_{\rho}}[\alpha
\cdot V(x) \cdot f ]- \frac{N^{3/2}}{2K+1}
            D_{\mu_\rho}^{K,x}(\sqrt{f})\}\label{SS2},
\end{align}
because of homogeneity of the system.  We next condition on the
density of particles in \reff{SS2} and apply the logarithmic
Sobolev inequality of Proposition \ref{lem:logsob} to obtain the
bound
\begin{align}
\Gamma^N_s \leq N^{3/2} \sup_{f}
\sup_{|\alpha|\leq\beta||J||_{\infty}}\sup_{|y|\leq
R}\{E_{\nu_{K,y}}[\alpha \cdot V(x) \cdot f ]- C
\frac{N^{3/2}}{K^3} H_{K,y}(f)\}\label{SS3},
\end{align}
for some constant $C$. Applying the entropy inequality again, we
obtain that for any $M$ positive the above is bounded by
\begin{eqnarray}
\hspace{-0.3cm}             N^{3/2} \sup_{f}
\sup_{|\alpha|\leq\beta||J||_{\infty}}\sup_{|y|\leq
R}\left\{\frac{1}{M}\log{E_{\nu_{K,y}}[\exp\{\alpha M V\}]}
            +\left(\frac{1}{M}- C \frac{N^{3/2}}{K^3}\right) H_{K,y}(f)\right\}.\label{line:lep3}
\end{eqnarray}
Choosing $M$ such that $\frac{1}{M}=C\frac{N^{3/2}}{K}$ and
letting $l=\frac{K}{\sqrt{N}}$, we may combine the last bound
above with \reff{line:lep1} to obtain
\begin{align}
E_{NEQ}[N^{1/2}|\int_0^tX_sds|] & \leq \frac{C t
        N^{1/2}}{\beta l^3} \sup_{|\alpha|\leq\beta}\sup_{|y|\leq R}
        \log{E_{\nu_{K,y}}[\exp\{\alpha
        Cl^3 V\}]}+\frac{C}{\beta}\label{line:lep5}.
\end{align}
The next lemma allows us to complete the argument.

\begin{lem}[Local Large Deviations]\label{LLD}

For any function $g$ which is at most linear in $k$, that is,
$g(k)\leq C (1+k)$, and for all $y\leq R$, there exists a constant
$C=C(R)$ such that for all $y$
\begin{eqnarray*}
\log E_{\nu_{K,y}}\left[\exp\left\{\gamma AV_{|x-y|\leq K}
\tilde{g}(\eta(x))\right\}\right] &\leq&
C\left\{\frac{\gamma^2}{K}+\frac{\gamma+\gamma^2}{K^2}\right\},\notag
\end{eqnarray*}
where $\tilde{g}=g(\eta(x))-\hat g(m^K(x))$.
\end{lem}
\bigskip
We apply this result to $E_{\nu_{K,y}}[\exp\{\alpha C l^3 V\}]$
and obtain that \reff{line:lep5} is bounded by
$C\{\beta^2l^2+\frac{l^{-2}+\beta^3l^4}{N^{1/2}}\}+\frac{C}{\beta}$,
from which the result follows.
\end{proof}

\begin{proof}[Proof of Lemma \ref{LLD}]

We begin by fixing a constant $M_0>0$.  For $R\geq y \geq M_0/K$,
we proceed as in \cite{Y92}. In \reff{line:LLD1}, \reff{line:LLD2}
and \reff{line:LLD3} below we look at the single site marginal of
the grand canonical measures. For simplicity of notation, we omit
the subscript $x$ from the functions $\eta$. We define $\bar{g}$
to be the centered version of $g$, $\bar{g}=g-E_{\mu_\rho}[g]$.

First, fix $\theta$ so that
$\int(\eta-\rho)e^{\theta(\eta-\rho)+\gamma/K\bar{g}(\eta)}d\mu_\rho(\eta)$
is equal to zero. A straightforward calculation shows that
\begin{align}
\theta=-\frac{\gamma}{K}\frac{<\eta;g>_{\rho}}{<\eta;\eta>_{\rho}}+O(\gamma^2/K^2).\label{line:LLD1}
\end{align}
As we are only considering values of $\rho$ inside a compact set,
these bounds are uniform. We use the notation $<g;f>_{\rho}$ to
denote the correlation of $g$ and $f$ under the measure
$\mu_{\rho}$.  Similarly, we define
$<f_1;f_2;f_3>_{\rho}=E_{\mu_{\rho}}[\bar{f_1}\bar{f_2}\bar{f_3}].$

Next we define the measure $P(\eta|\theta, g, \rho)$ by
\begin{eqnarray*}
&dP(\eta)= Z^{-1}(\theta, g, \rho)
\exp\{\theta(\eta-\rho)+\gamma/K \bar g(\eta) \}d\mu_\rho(\eta),&
\end{eqnarray*}
with normalizing factor, $Z(\theta, g, \rho)$, given by
\begin{eqnarray}
\int e^{\theta(\eta-\rho)+\gamma/K\bar{g}(\eta)}d\mu_\rho(\eta)=
1+\gamma^2/K^2 H_2+O(\gamma^3/K^3)\label{line:LLD2},
\end{eqnarray}
where
$H_2=-\frac{<\eta,g>_{\rho}^2}{<\eta;\eta>_{\rho}}+<g;g>_{\rho}$.
\begin{rem}\label{rem:lepfem}
Note that it is because of this step that we restrict ourselves to
linear functions $g$, as otherwise the normalizing function will
not be finite.
\end{rem}
We also define the variance under our modified measure:
\begin{eqnarray}
\sigma^2(\theta,g,\rho)&=&Z^{-1}(\theta,g,\rho)\int(\eta-\rho)^2e^{\theta(\eta-\rho)+\gamma/K\bar{g}(\eta)}d\mu_\rho(\eta)\notag\\
&=&\sigma^2\,Z^{-1}(\theta,g,\rho)\left\{1-2\gamma/KH_1+O(\gamma^2/K^2)\right\},\label{line:LLD3}
\end{eqnarray}
where $\sigma^2=<\eta,\eta>_\rho$, and
$H_1=\frac{<\eta;g>_{\rho}<\eta;\eta;\eta>_{\rho}}{<\eta;\eta>_{\rho}^2}-\,\frac{<\eta;\eta;g>_{\rho}}{<\eta;\eta>_{\rho}}$.
We now re-write the quantity of interest using the notation
developed above.
\begin{eqnarray}
\hspace{-0.5cm}E_{\nu_{K,y}}[\exp\{\gamma AV_{|x-y|\leq K} \bar{g}(\eta(x))\}]&=&\frac{\int_{\bar{\eta}=\rho}e^{\theta\sum_{|x-y|\leq K}(\eta-\rho)+\gamma/K\bar{g}(\eta)}d\mu_\rho(\eta)}{\int_{\bar{\eta}=\rho}d\mu_\rho(\eta)}\notag\\
&=&Z^{-K}(\theta,g,\rho)\cdot\frac{P(\bar{\eta}=\rho|\theta,g,\rho)}{P(\bar{\eta}=\rho|\theta=g=0,\rho)}.\label{line:LLD4}
\end{eqnarray}
We now apply Proposition \ref{thm:localCLT} with $J=8$ to show
that
\begin{eqnarray*}
\frac{\sqrt{K\sigma^2(\theta,g,\rho)}P(\bar{\eta}=\rho|\theta,g,\rho)}{\sqrt{K\sigma^2(\rho)}P(\bar{\eta}=\rho|\theta=g=0,\rho)}&=&1+O(\gamma/K^2).
\end{eqnarray*}
We combine this last bound with \reff{line:LLD1} through
\reff{line:LLD3} to show that \reff{line:LLD4} behaves like
$1+\gamma/KH_1-1/2\gamma^2/KH_2+O((\gamma+\gamma^2)/K^2).$

It remains to study the effect of considering $\tilde{g}$ vs.
$\bar{g}$ in the expectation under study.  Let $g^* =
AV_{|x-y|\leq K}\bar{\tilde{g}}(\eta(x))$.  By Jensen's inequality
we have for any positive $\alpha$:
\begin{eqnarray*}
-\alpha^{-1}\log E_{\nu_{K,y}}[\exp\{-\alpha g^*\}] \leq
E_{\nu_{K,y}}[g^*] \leq \alpha^{-1}\log E_{\nu_{K,y}}[\exp\{\alpha
g^*\}].
\end{eqnarray*}
From previous bounds we have
\begin{eqnarray*}
\hspace{-0cm}\alpha^{-1}\log E_{\nu_{K,y}}[\exp\{\pm \alpha
AV_{|x-y|\leq K}\bar{\tilde{g}}(\eta(x))\}] =\pm
1/KH_1-1/2\alpha/KH_2+O((1+\alpha)/K^2),
\end{eqnarray*}
which equals $H_1/K+O(K^2)$ if we select $\alpha=1/K$.  
We combine these results to obtain the required result for the
case $y\geq\frac{M_0}{K}$.
\begin{eqnarray*}
\hspace{-0cm}\gamma^{-1}\log\left(E_{\nu_{K,y}}\left[\exp\left\{\gamma
AV_{|x-y|\leq K}\tilde{g}(\eta(x))\right\}\right]\right)
=\gamma/KH_2+O\left((1+\gamma)/K^2\right).
\end{eqnarray*}
For $ y < M_0/K$ the bound is much simpler as in this case the
total number of particles, and hence any \emph{non-local}
function, is bounded.  In this case the result follows by a simple
Taylor expansion.
\end{proof}

\subsubsection{Equivalence of Solutions.}

By the hydrodynamic scaling limit we know that $m^K$ and $\rho$
solve the same limiting differential equations.  We use here a
standard PDE trick to show that the solutions themselves are
close.  This idea was first developed in \cite{Y92}.

\begin{lem}\label{lem:hnorm}
Under the assumptions of this chapter we have that
\begin{align*}
\lim_{l\rightarrow0}\lim_{N\rightarrow\infty}E[\int_0^t
N^{-1/2}\sum_{x\in \bb T^N}(m^K(x)-\rho(x))^2ds]=0.
\end{align*}
\end{lem}

\begin{proof}
Let $(\triangle f)(x) = f(x+1)-2f(x)+f(x-1)$ denote the discrete
Laplacian, and write its inverse as the matrix $G=\{g(x,y)\}$.  As
we work with centered functions on the discrete torus, this
inverse is well-defined.

A careful calculation gives that $L||\bar{Y}^N_t||_{*}^2$ is equal
to
\begin{align}
&L \left\{-\frac{1}{N}\sum_{x\in \bb T^N}
        \sum_{y\in \bb T^N}N^{-2}g(x,y)\bar{Y}^N_t(x)\tilde{Y}^N_t(y)\right\}\notag\\
& = -\sqrt{N}<c(\eta^N)-\varphi(\rho),\bar{Y^N}> + \sum_{x\in \bb
T^N}
    c(\eta^N(x))+\frac{1}{N}\sum_{x\in \bb T^N}c(\eta^N(x)).\label{line:hnorm1}
\end{align}
Similarly, we obtain
\begin{align}
L<\bar{\eta}^N-\bar{\rho},\bar{Y}^N_t>_{*}
    & =
    \frac{\sqrt{N}}{2}<\bar{\eta^N}-\bar{\rho},c(\eta^N)-\varphi(\rho)>.\label{line:hnorm2}
\end{align}
Combining \reff{line:hnorm1} with \reff{line:hnorm2} we get
\begin{align}
&E_{NEQ}[N^{-1/2}||\bar{Y}^N_t||_{*}^2+2<\bar{\eta}^N-\bar{\rho},\bar{Y}^N_t>_{*}]\notag\\
&\hspace{1cm}- E_{NEQ}[N^{-1/2}||\bar{Y}^N_0||_{*}^2+2<\bar{\eta}^N-\bar{\rho},\bar{Y}^N_0>_{*}]\notag\\
  &= \hspace{1cm}E_{NEQ}\left[\int_0^t N^{-1/2}H_s(x)ds\right]+\,\,\,N^{-1/2}E_{NEQ}\left[\int_0^t AV_{x\in \bb T^N}
  c(\eta_s^N(x))ds\right]\label{line:hnorm3}
\end{align}
where we let
$H_s(x)=c(\eta_s^N(x))-\left\{c(\eta^N_s(x))-\varphi(\rho_s(x))\}\{\eta_s^N(x)-\rho_s(x)\right\}.$

We next replace $H(x)$ with its local average
$$\bar{H}^K(x)=AV_{|x-y|\leq
K/N}\left[\rule[-0.3cm]{0cm}{0.7cm}\,c(\eta(y))-\{c(\eta(y))-\varphi(\rho(x))\}\{\eta(y)-\rho(x)\}\right].$$
Assumptions (LG) and (M) imply that we have finite exponential
moments for the measure $\mu_\rho$.  By the entropy inequality,
and assumption (H1) we conclude that
\begin{eqnarray}\label{SSS}
N^{-1/2}E_{NEQ}[\int_0^tAV_{x\in \bb
T^N}c(\eta^N_s(x))ds]=O(N^{-1/2}).
\end{eqnarray}
Thus, using summation by parts, along with $l=\frac{K}{\sqrt{N}}$,
we obtain that
\begin{align*}
E_{NEQ}\left[\int_0^t N^{-1/2}\sum_{x\in \bb T^N}
(H_s(x)-\bar{H}_s^K(x))(s)ds\right] = O(l).
\end{align*}
Combining the above result with \reff{line:hnorm3} we have
\begin{align}
&E_{NEQ}[N^{-1/2}||\bar{Y}^N_t||_{*}^2+2<\bar{\eta}-\bar{\rho},\bar{Y}^N_t>_{*}] - E[N^{-1/2}||\bar{Y}^N_0||_{*}^2+2<\bar{\eta}-\bar{\rho},\bar{Y}^N_0>_{*}]\notag\\
&= E_{NEQ}\left[\int_0^t N^{-1/2}\sum_{x\in \bb
T^N}\bar{H}^K_s(x)ds\right] +
  O(l)+O(N^{-1/2}).\label{line:hnorm4}
\end{align}
We will make use of the local equilibrium principle to handle the
remaining term.  Define $F(\eta(y),\rho(x))$ to be the quantity
$$c(\eta(y))-\{c(\eta(y))-\varphi(\rho(x))\}\{\eta(y)-\rho(x)\},$$
and let $\hat{F}(m^K(x),\rho(x))$ denote its expectation under the
measure $\mu_{\rho=m^K(x)}$.  A brief calculation shows that $\hat
F$ is equal to
$\left\{\rule[-0.3cm]{0cm}{0.7cm}\varphi(m^K(x))-\varphi(\rho(x))\right\}\left\{\rule[-0.3cm]{0cm}{0.7cm}m^K(x)-\rho(x))\right\}.$
Lastly, let $s^K(x)=AV_{|x-y|\leq \frac{K}{N}}\,\eta^2(x)$.  We
split $\bar{H}^K(x)$ into three pieces:
\begin{eqnarray*}
\bar{H}^K(x)&=&AV_{|x-y|\leq \frac{K}{N}}\left[\rule[-0.3cm]{0cm}{0.7cm}F(\eta(y),\rho(x))\right.\\
&&\hspace{1.2in}\left.\rule[-0.3cm]{0cm}{0.7cm}-\hat{F}(m^K(x),\rho(x))\right]\bb{I}(m^K(x)\leq R \cap s^K(x)\leq R_2)\\
&&+\,\,AV_{|x-y|\leq \frac{K}{N}}\,\,F(\eta(y),\rho(x)) \,\,\bb{I}(m^K(x) > R \cup s^K(x)> R_2)\\
&&+\,\,AV_{|x-y|\leq
\frac{K}{N}}\,\,\hat{F}(m^K(x),\rho(x))\,\, \bb{I}(m^K(x)\leq R \cap s^K(x)\leq R_2)\\
&=&AV_{|x-y|\leq
\frac{K}{N}}\left[\rule[-0.3cm]{0cm}{0.7cm}F^{(1)}(x)+F^{(2)}(x)+F^{(3)}(x)\right].
\end{eqnarray*}
 By the local
equilibrium principle the term involving $F^{(1)}$ vanishes in the
limit. Even though the term $F-\hat{F}$ is not linear in $\eta$ as
required in the local equilibrium principle, the addition of the
condition $\bb{I}(s^K(x)\leq R_2)$ fixes the problem pointed out
in Remark~\ref{rem:lepfem} . 
Applying bounds from Remark~\ref{rem:linerate} of Section
\ref{ch:def} we have that
\begin{align}
AV_{|y-x|\leq \frac{K}{N}} F(\eta(y),\rho(x))
\leq\tilde{c}\,\,m^K(x)-c_1AV_{|y-x|\leq
\frac{K}{N}}(\eta(y)-\rho(x))^2,\label{line:es11}
\end{align}
for some constant $\tilde{c}>0$. If $m^K(x)>R$, we fix $C$ and
choose $R$ so that $2||\rho||_{\infty}<R-C$, which implies that
$m^K(x)-\rho(x)>C$ and hence \reff{line:es11} is smaller than
\begin{eqnarray*}
-A\,(m(x)-\rho(x))^2 \bb{I}(m^K(x)>R)
\end{eqnarray*}
We next pick $C$ large enough so that $A=c_1-2\tilde{c}/C>0$.

If, on the other hand, $m^K(x)\leq R$, then we control
\reff{line:es11} using $s^K(x)>R_2$. For $R_2$ large enough
\begin{align*}
AV_{|y-x|\leq K/N}(\eta(y)-\rho(x))^2>R_2-2R||\rho||_{\infty}.
\end{align*}
Hence, for $c>0$ take $R_2>R/c+2R||\rho||_{\infty}$, we get that
$R<c\,(R_2-2||\rho||_{\infty}R)$ implying
$$\tilde{c}\,m^K(x)\leq \tilde{c}\,R \leq \tilde{c}c\, AV_{|y-x|\leq K/N}(\eta(y)-\rho(x))^2.$$
We thus obtain a bound on \reff{line:es11} of
\begin{eqnarray*}
-B (m^K(x)-\rho(x))^2\bb{I}(m^K(x)\leq
R,m^K_2(x)>R_2)
\end{eqnarray*}
where $B=c_1-\tilde c c$ is positive as long as $c$ was chosen to
be sufficiently small. Thus, we obtain that for fixed $R$ and
$R_2$ large enough, and $\tilde{A}=min(A,B)>0$
\begin{align}
AV_{|x-y|\leq K/N}F^{(2)}(x)\leq
-\tilde{A}(m^K(x)-\rho(x))^2\bb{I}(m^K(x)\leq
R,m^K_2(x)>R_2)\label{line:hnorm6}
\end{align}
A similar bound holds for $F^{(3)}$, using the fact that
$\varphi'\geq c_1$. Combining these last bounds together with
\reff{line:hnorm4} we find
\begin{eqnarray}
&&\hspace{-1cm}E_{NEQ}[N^{-1/2}||\bar{Y}^N_t||_{*}^2+2<\bar{\eta}-\bar{\rho},\bar{Y}^N_t>_{*}] - E[N^{-1/2}||\bar{Y}^N_0||_{*}^2+2<\bar{\eta}-\bar{\rho},\bar{Y}^N_0>_{*}]\notag\\
  &\leq&-A^*E_{NEQ}[\int_0^t N^{-1/2}\sum_{x\in \bb T^N}(m^K(x)-\rho(x))^2 \bb{I}(m^K(x)>R\cup m^K_2(x)>R_2)ds]+o(1),\notag
\end{eqnarray}
where $A^* = \tilde A +c_1$.  We shall now use assumption (F1)
together with the inequalities
\begin{eqnarray}\label{eos:line}
\pm 2 <\bar{\eta}-\bar{\rho},\bar{Y}_t^N>_{*}\leq
2N^{1/2}||\bar{\eta}-\bar{\rho}||^2_{*}
+\frac{1}{2}N^{-1/2}||\bar{Y}_t^N||^2_{*}
\end{eqnarray} and
$||\bar{\eta}-\bar{\rho}||^2_{*}\leq C
||\bar{\eta}-\bar{\rho}||^2,$ to conclude that
\begin{align*}
\frac{1}{2}N^{-1/2}E_{NEQ}[||\bar{Y}^N_t||_{*}^2]
\leq-c\,E_{NEQ}[\int_0^t N^{-1/2}\sum_{x\in \bb
T^N}(m^K(x)-\rho(x))^2 ds]+o(1).
\end{align*}
The result follows.
\end{proof}

\subsubsection{Conclusion.}

By Lemmas \ref{lem:LEP} and \ref{lem:hnorm} we have  that
\reff{BGtodo} holds for $i=2,3$. We consider $\Phi_1^N$ next. By
summation by parts and \reff{SSS} we have for $f\in C^1(\bb T)$
\begin{eqnarray}\label{BGbound1}
|E[\int_0^t<f,\Phi_1^N(s)>ds]|\leq C l,
\end{eqnarray}
Using arguments similar to those used to obtain \reff{line:hnorm6}
we show that for $i=4,5$
\begin{eqnarray}\label{BGbound45}
|<1,\Phi^N_i(x)>| \leq C' (m^K(x)-\rho(x))^2\bb{I}[m^K_t(x)> R]
\end{eqnarray}
for some positive constant $C'$ and then apply
Lemma~\ref{lem:hnorm}. This concludes the proof of
Theorem~\ref{BG}.

\subsection{Quadratic Variation}\label{section:QV}
We need to show that $Y_t$ has quadratic variation given by
\reff{limitqv1}. We know that under  $P_N$, $N^N_t(f)$ defined in
\reff{line:qvcalc} is a martingale.  Hence, it remains to show the
following proposition. Note that it implies that test function $f$
must be at least $C^2(\bb T)$.
\begin{prop}
For a function $h$ in $C^1(\bb T)$.
\begin{align}
E\left[\left|\int_0^t<h,c(\eta^N_s)-\varphi(\rho_s)>ds\right|\right]\rightarrow
0 .
\end{align}
\end{prop}

\begin{proof}
This follows from the Boltzmann-Gibbs principle of Theorem
\ref{BG}, which is a stronger result, in combination with the
hydrodynamic scaling limit of Proposition~\ref{HSL}.
\end{proof}

\subsection{Unique Solution to the Martingale Problem.\label{sec:uniqueone}}

We recall here a result due to Holley and Stroock \cite{Ho-S}
which guarantees the existence of a weakly unique solution to our
martingale problem.

Let $\mathcal{B}_t$ denote the nonpositive operator $\frac{1}{2} D
(\rho_t)\triangle$ and $\mathcal{P}^D(t,s)$ the associated
semigroup.  Also, let $\mathcal{S}_t$ be the operator
$\sqrt{\varphi(\rho_t)}\nabla$.

\begin{thm}\label{uniqueone}
Fix a positive integer $m\geq2$.  Let $P$ be a probability measure
on the space $\{(C[0,T],H_{-m}),\mathcal{ F}\}$, where
$\mathcal{F}=\cup_{t\geq0}\mathcal{F}_t$, and $\mathcal{F}_t$ is
the canonical filtration of the process $<Y_t,f>$. Assume that for
any $f \in C^{\infty}(\mathbb{T})$
\begin{eqnarray}
\begin{array}{cl}
&<M_t,f>\,=\,<Y_t,f>-<Y_0,f>-\int_0^t<Y_s,\mathcal{B}_sf>ds\\
\mbox{and}&<M_t,f>^2-\int_0^t||\mathcal{S}_sf||^2_2ds
\end{array}
\label{condmart}
\end{eqnarray}
are $L^1(P)$ martingales with respect to $\mathcal{F}_t$.  Then,
for all $0\leq s\leq t$, $f$ in $C^{\infty}(\mathbb{T})$, and
subsets $A$ of $\mathbb{R}$,
\begin{eqnarray*}
&&\hspace{-2cm}P[<Y_t,f>\in A|\mathcal{F}_s]\\
&= &\int_A \frac{1}{\sqrt{2\pi \int_s^t
||\mathcal{S}_u\mathcal{P}^D(t,u)f||^2_2du}}\exp\left\{-\frac{(y-<Y_s,\mathcal{P}^D(t,s)f>)^2}{2\int_s^t
||\mathcal{S}_u\mathcal{P}^D(t,u)f||^2_2du}\right\}dy,
\end{eqnarray*}
$P$ almost surely.  In particular, \reff{condmart} together with a
unique initial distribution for $<Y_0,f>$ uniquely determines $P$
on $\{(C[0,T],H_{-m}),\mathcal{ F}\}$.
\end{thm}

Thus, to guarantee uniqueness of the limiting measure $P$, we need
only check that we have a unique initial measure; this is simply
condition (F2).

\subsection{Tightness\label{sec:tightone}}

Let $P_N$ denote the probability measure on $D([0,T],
\mathcal{H}_{-m})$ induced by the fluctuation field $Y_t^N$
started in the measure $\zeta_0^N$, and $E_{NEQ}$ the associated
expectation. In this section we prove that the sequence of
probability measures $P_N$ is tight. A consequence of the proof
will be that the limit points are concentrated on the space of
continuous paths $C([0,T], \mathcal{H}_{-m})$.

For a function $F$ in $D([0,T], \mathcal{H}_{-m})$, define the
(uniform) modulus of continuity, $\omega_\delta(F)$, for a fixed
$\delta>0$:
\begin{eqnarray*}
\omega_\delta(F)=\sup_{|s-t|\leq\delta, 0\leq s,t \leq
T}||F_t-F_s||_{-m}.
\end{eqnarray*}
To simplify notation we will often denote this supremum simply as
$\sup_{s,t}||F_t-~F_s||_{-m}.$ By the well-known Arzela-Ascoli
result, it follows that to show the sequence $P_N$ is tight we
need to show that it satisfies two conditions:
\begin{enumerate}
\item[](T1).\hspace{1cm}$\lim_{M \rightarrow \infty}\limsup_{N \rightarrow \infty}
P_N [\sup_{0\leq t \leq T} ||Y_t||_{-m}>M] =0$
\item[] (T2).\hspace{1cm}$\lim_{\delta \rightarrow 0}\limsup_{N \rightarrow
\infty}P_N[\omega_\delta(Y)> \epsilon]=0 \,\,\,\,\,\,\,\,\,\,\,\,
\forall \epsilon
>0.$
\end{enumerate}

We will also make use of the result, due to Aldous, \cite{Ald}.
For a proof, see for example Proposition 1.6, of Section 4  in
\cite{K-L}. We will state it in some generality.  Let $\mc X$
denote a Polish space with metric $d(\cdot, \cdot)$.  For a
function $g$ in $D([0,T],\mc  X)$, define the modified modulus of
continuity, $\omega'_\delta(g)$, for a fixed $\delta>0$:
\begin{eqnarray*}
\omega'_\delta(g)=\inf_{\{t_i\}_{i=0}^r} \max_{0\leq i <
r}\sup_{t_i\leq s<t<t_{i+1}}d(g_t, g_s),
\end{eqnarray*}
where the infimum is taken over all partitions $\{t_i\}_{i=0}^r$
of $[0,T]$ such \linebreak that $t_0=0 < t_1 < \hdots < t_r=T$ and
$t_{i+1}-t_i>\delta$.

\begin{prop}\label{Aldous}
A sequence of probability measures $P_N$ on $D([0,T], \mc X)$
satisfies the condition
$$\lim_{\delta \rightarrow 0}\limsup_{N
\rightarrow \infty}P_N[\omega'_\delta(Y)> \epsilon]=0$$ provided
that
\begin{eqnarray}\label{aldouscondition}
\lim_{\delta \rightarrow 0}\limsup_{N \rightarrow
\infty}\sup_{\stackrel{\tau \in \mc T}{\theta \leq \delta}}
P_N[d(Y_\tau,Y_{\tau + \theta})> \epsilon]=0,
\end{eqnarray}
where $\mc T$ denotes the collection of all stopping times bounded
by $T$.
\end{prop}
\begin{prop}\label{prop:tightone}
The sequence of measures $P_N$ is tight in $D([0,T],
\mathcal{H}_{-m})$, for any $m\geq4$.  Moreover, all limit points
are concentrated on continuous paths.
\end{prop}
We first consider condition (T2). Recall that $\tilde{Y}_s^N$ is
the field
 $\sqrt{N}(c(\eta^N_t(x)-\rho_t(x)).$ We may hence write
\begin{eqnarray*}
||Y^N_t-Y^N_s||_{-m}&\leq&\sup_{||f||_{m}\leq 1}|<\triangle_Nf,
\int_s^t \tilde{Y}_u^Ndu
>|+||M^N_t-M^N_s||_{-m}
\end{eqnarray*}
where $<M_t^N,f>$ is a martingale with quadratic variation given
by the martingale in~\reff{line:qvcalc}. Using this result, we
reduce the proof of condition (T2) to showing that both quantities
\begin{enumerate}
\item[]$(T2)_A$ \hspace{1cm} $P_N\left[\sup_{s,t}\sup_{||f||_{m}\leq 1}|<\triangle_Nf, \int_s^t \tilde{Y}_u^Ndu
>|> \frac{\epsilon}{2}\right]$
\item[]$(T2)_B$\hspace{1cm} $P_N\left[\sup_{s,t} ||M^N_t-M^N_s||_{-m} > \frac{\epsilon}{2}\right]$.
\end{enumerate}
decrease as $|t-s|\leq\delta$ converges to zero.  We begin with
$(T2)_B$, and we split the proof into several lemmas.

\begin{lem}[Bounding the tails 1.]\label{TBboundtails}
There exists a finite constant $C$ such that for every
eigenfunction $e_z$ of $\mc H_m$,
$$\limsup_{N\rightarrow\infty}E_{NEQ}\left[\sup_{0\leq t\leq T}|<M_t^N,e_z>|^2\right]\leq C T \{1+<\triangle e_z,\triangle e_z>\}$$
\end{lem}

\begin{lem}[Bounding the tails 2.]\label{TBboundtails2}
For $m\geq 3$,
$$\lim_{K\rightarrow\infty}\limsup_{N\rightarrow\infty}E_{NEQ}\left[\sup_{0\leq t\leq T}\sum_{z\geq K}\{<M_t^N,e_z>\}^2\gamma^{-m}_z\right]=0.$$
\end{lem}
This last lemma implies that to prove condition $(T2)_B$ is
satisfied, we need only check that for bounded $K$
$$\lim_{\delta\rightarrow
0}\limsup_{N\rightarrow\infty}P_N\left[\sup_{s,t}\sum_{|z|\leq
K}\{<M_t^N-M_s^N,e_z>\}^2\gamma^{-m}_z>\epsilon'\right]=0.$$ This
will follow if we can prove the following result.
\begin{lem}\label{martboundsone} For any $f$ in $C^3(\bb T)$
$$\lim_{\delta\rightarrow
0}\limsup_{N\rightarrow\infty}P_N\left[\sup_{s,t}|<M_t^N-M_s^N,f>|>\epsilon'\right]=0.$$
\end{lem}

We make use of the following fact
\begin{eqnarray*}
\omega_{\delta} (<M_t^N,f>)\leq 2 \omega'_{\delta}(<M_t^N,f>)+
\sup_{t}|<M_t^N,f>-<M_{t^-}^N,f>|.
\end{eqnarray*}
By the definition of $M_t^N$ we know that
$$|<M_t^N,f>-<M_{t^-}^N,f>|=|<Y_t^N,f>-<Y_{t^-}^N,f>|,$$ and notice
that because at most one particle makes a jump at any given time,
this quantity is bounded above by $||f'||_{\infty}N^{-3/2}.$  We
therefore obtain that
\begin{eqnarray}\label{aldouscondition2}
\omega_{\delta} (<M_t^N,f>)\leq 2 \omega'_{\delta}(<M_t^N,f>)+
||f'||_{\infty}N^{-3/2}
\end{eqnarray}

\begin{proof}[Proof of Lemma \ref{martboundsone}.]
By \reff{aldouscondition2} above and Proposition \ref{Aldous} it
is enough to check that
\begin{eqnarray}
\lim_{\delta \rightarrow 0}\limsup_{N \rightarrow
\infty}\sup_{\stackrel{\tau \in \mc T}{\theta \leq \delta}}
P_N\left[|<M_{\tau +\theta}^N-M_{\tau}^N,f>|>\epsilon'\right]=0,
\end{eqnarray}
By Chebychev's inequality, and because both $\tau$ and $\theta$
are bounded stopping times, we may bound $P_N[|<M_{\tau +\theta
}^N-M_{\tau}^N,f>|>\epsilon']$ by
$$\frac{1}{2\epsilon^2}E_{NEQ}\left[\int_0^\delta \sum_{x\in
       \bb{T}^N}[c(\eta^N_s(x))][\nabla_Nf(x)^2+\nabla_Nf(x+1/N)^2]ds\right].$$
Applying the results of Section \ref{section:QV} completes the
proof.
\end{proof}

We finish the argument for $(T2)_B$ by proving Lemmas
\ref{TBboundtails} and \ref{TBboundtails2}.

\begin{proof}[Proof of Lemma \ref{TBboundtails}]
Because $M_t^N$ is a martingale we may use Doob's inequality to
bound $E_{NEQ}[\sup_{0\leq t\leq T}|<M_t^N,e_z>|^2]$ by
$4E_{NEQ}[|<M_T^N,e_z>|^2] $ and this is equal to
$$\frac{4}{2}\int_0^T
\sum_i\sum_{x\in
       \bb{T}^N}[c(\eta^N_s(x))][\nabla_Ne_z(x)^2+\nabla_Ne_z(x-1/N)^2+\nabla_Ne_z(x+1/N)^2]ds.$$
Using the convergence established in section \ref{section:QV}, and
the fact that $\varphi$ is bounded uniformly, we may bound the
above by $T ||\varphi||_\infty
\{1+<\triangle^Ne_z,\triangle^Ne_z>\} $.
\end{proof}

\begin{proof}[Proof of Lemma \ref{TBboundtails2}]
This follows since by Lemma \reff{TBboundtails} the quantity
$$E_{NEQ}\left[\sup_{0\leq t\leq T}\sum_{z\geq K}\{<M_t^N,e_z>\}^2\gamma^{-m}_z\right]$$
is bounded by $C \sum_{|z|\geq K}
\gamma_z^{-m}\{1+<\triangle^Ne_z,\triangle^Ne_z>\}.$
\end{proof}

We next turn to condition $(T2)_A$:
$$(T2)_A \hspace{1cm}\lim_{\delta\rightarrow 0}\limsup_{N\rightarrow\infty} P_N\left[\sup_{s,t}\sup_{||f||_{m}\leq 1}\left|\int_s^t <\tilde{Y}_u^N,\triangle_Nf>du
\right|> \frac{\epsilon}{2}\right]=0,$$ where $\tilde{Y}_s^N$ is
the field $\sqrt{N}(c(\eta^N_t(x)-\varphi(\rho_t(x))).$  Similarly
to the proof of the Boltzmann-Gibbs principle, we split this field
up as the sum of the following terms:
\begin{eqnarray}
\left.\begin{array}{ccl}
\Psi^N_{1,t}&=&\sqrt{N}\{c(\eta^N_t(x))-\bar{c}_K(\eta^N_t(x))\}\\
\Psi^N_{2,t}&=&\sqrt{N}\{\bar{c}_K(\eta^N_t(x))-\varphi(m^K_t(x))\}\bb{I}[m^K_t(x)\leq R]\\
\Psi^N_{3,t}&=&\sqrt{N}\varphi'(\rho_t(x))\{m^K_t(x)-\rho_t(x)\}\bb{I}[m^K_t(x)\leq R]\\
\Psi^N_{4,t}&=&\sqrt{N}\{\varphi(m^K_t(x))-\varphi(\rho_t(x))-\varphi'(\rho_t(x))\{m^K_t(x)-\rho_t(x)\}\}\bb{I}[m^K_t(x)\leq R]\\
\Psi^N_{5,t}&=&-\sqrt{N}\varphi(\rho_t(x))\bb{I}[m^K_t(x)> R]\\
\Psi^N_{6,t}&=&\sqrt{N}\bar{c}_K(x)\bb{I}[m^K_t(x)>
R]\end{array}\right.\label{chooseK}
\end{eqnarray}
It hence remains to prove that for $i=1,\ldots,6$
$$\lim_{\delta\rightarrow 0}\limsup_{N\rightarrow\infty} P_N\left[\sup_{s,t}\sup_{||f||_{m}\leq 1}\left|\int_s^t <\Psi_{i,u}^N,\triangle_Nf>du
\right|> \frac{\epsilon}{12}\right]=0$$ For $i\neq2,3$ we first
apply Chebychev's inequality and use Lemma \ref{Aldous}.  For
$i=1$ we argue as in \reff{BGbound1}.  For $i=4$ we have Lemma
\ref{lem:hnorm}. Lastly, we use \reff{BGbound45} and Lemma
\ref{lem:hnorm} to handle $i=5,6$.
\begin{eqnarray*}
P_N\left[\sup_{||f||_{m}\leq 1}\left|\int_\tau^{\tau+\delta}
<\Psi_{i,u}^N,\triangle_Nf>du \right|> \epsilon'\right]\leq
\frac{1}{\epsilon'}E_{NEQ}\left[\sup_{||f||_{m}\leq 1}\int_0^T
\left|<\Psi_{i,u}^N,\triangle_Nf>\right|du \right].
\end{eqnarray*}
It remains to handle the terms with $i=2,3$. In order to do this
we use a special case of the well-known Garcia-Rodemich-Rumsey
inequality.
\begin{lem}[Garcia-Rodemich-Rumsey inequality]\label{GRR} Given
a function $g$ and a strictly increasing function $\psi$ such that
$\psi(0)=0$ and $\lim_{u \rightarrow\infty}\psi(u)=\infty$, then
\begin{align*}
|g(t)-g(s)|\leq 8
\int_0^\delta\psi^{-1}(\frac{4B}{u^2})\frac{du}{\sqrt{u}},\\
\mbox{where}\hspace{1cm}B=\int_0^T\int_0^T
\psi(\frac{|g(t)-g(s)|}{\sqrt{t-s}})dtds.
\end{align*}
\end{lem}
We choose $\psi(x)=e^{N|x|}-1$ in the above inequality.  By
integration by parts we obtain
\begin{align*}
&\int_0^\delta \log  \left(1+\frac{4B}{u^2}\right)
\frac{1}{\sqrt{u}}du \leq \sqrt{\delta}\{\log
(1+\frac{4B}{\delta^2})+4\}.
\end{align*}
To simplify notation we drop the subscripts on the field
$\Psi^N_{i,s}$ and simply write it as $\Psi(s)$.  We apply the
Garcia-Rodemich-Rumsey inequality to any such field $\Psi(s)$ in
the following way: choose $g(t)=\int_0^t <\Psi(u),f> du$, where
$||f||_{m}\leq 1$. We then obtain
\begin{align*}
|\int_s^t <\Psi(u),f> du| \leq  \sqrt{\delta}\{\log
(1+\frac{4B}{\delta^2})+4\},
\end{align*}
and we may also bound $B$ by
\begin{align*}
\int_0^T \int_0^T \left\{\exp\left( N \frac{||\int_s^t
\Psi(u)du||_{-m}}{\sqrt{t-s}}\right) -1\right\} ds dt,
\end{align*} and we have the inequality
\begin{align}\label{defA0}
A_0=\sup_{s,t}||\int_s^t \Psi(u)du ||_{-m} \leq
\frac{8}{N}\delta^{-1/2}\{\log(1+4B_0\delta^{-2})+4\},
\end{align}
where
\begin{eqnarray}\label{defB0}
B_0=\int_0^T \int_0^T \exp\left( N \frac{||\int_s^t
\Psi(u)du||_{-m}}{\sqrt{t-s}}\right)ds dt-T^2.
\end{eqnarray}
It follows that
\begin{eqnarray*}
B_0\leq \frac{\delta^2}{4}e^{N\epsilon /16 \sqrt{\delta}} \,\,
\Rightarrow A_0 \leq \epsilon /2 + \frac{64}{N},
\end{eqnarray*}
and this last quantity is smaller than $\epsilon$ for large enough
$N$. We thus have that for any field $\Psi$ the following
inequality holds
\begin{eqnarray}
P_N(A_0 > \epsilon ) \leq P_N\left(B_0 >
\frac{\delta^2}{4}e^{N\epsilon /16
\sqrt{\delta}}\right).\label{ineqBA1}
\end{eqnarray}
Next, we apply the entropy inequality (see for example
Prop. 8.2, of Appendix A in \cite{K-L}) to bring the problem back
into equilibrium measure. We also use here assumption (H1).  We
thus obtain that for any set $A$ such that $P_{EQ}^N(A)\leq e^{-N
C(\delta)}$
\begin{eqnarray}
P_N[A]&\leq& [\log 2 + H_N(0) ][\log(1+P_{EQ}^N(A)^{-1})]^{-1},\notag\\
      &\leq& [CN][\log(1+P_{EQ}^N(A)^{-1})]^{-1},\notag\\
      &\leq& 1/C(\delta)\label{ineqBA2}
\end{eqnarray}
where $P_{EQ}^N$ is the measure induced by the field $Y^N$ started
in the equilibrium measure.  We have thus reduced the problem to
showing
\begin{eqnarray}\label{tightone:line1}
P_{EQ}^N\left(B_0 > \frac{\delta^2}{4}e^{N\frac{\epsilon}{16}
\sqrt{\delta}}\right) \leq e ^{-N C(\delta)},
\end{eqnarray}
where $C(\delta)\rightarrow\infty$ as $\delta \rightarrow 0$. We
now apply Chebychev's inequality
\begin{eqnarray*}
P_{EQ}^N\left(B_0>\frac{\delta^2}{4}\exp{\left\{N\frac{\epsilon}{16}\delta^{-1/2}\right\}}\right)\leq
4 \frac{E_{EQ}[B_0]}{\delta^2}\exp\left\{-N\frac{\epsilon}{16}
\delta^{-1/2}\right\}.
\end{eqnarray*}
Suppose next that we have that
\begin{align}\label{tightabove}
E_{EQ}[B_0]\leq e^{C N},
\end{align}
where $C$ is some positive constant.  Notice that $x \leq e^{-x}$
for small $x$, and hence if \reff{tightabove} holds, we have that
there exists a $\delta_0>0$ and a positive constant $C'$ depending
on $\delta_0$ such that for all $\delta<\delta_0$ we have
\begin{eqnarray}
P_{EQ}^N\left(B_0>\frac{\delta^2}{4}\exp{\left\{N\frac{\epsilon}{16}\delta^{-1/2}\right\}}\right)\leq
\exp\left\{-N C' \frac{1}{\delta^2}-
N\frac{\epsilon}{16}\delta^{-1/2} \right\}.
\end{eqnarray}
Hence \reff{tightone:line1} holds.  We have in fact proved the
following.

\begin{prop} \label{usefultight}
Suppose that there exists a constant $C>0$ and an $N_0$ such that
for all $N>N_0$
\begin{align}\label{tightone:line3}
E_{EQ}[\int_0^T \int_0^T \exp\left( N \frac{||\int_s^t
\Psi(u)du||_{-m}}{\sqrt{t-s}}\right)ds dt]\leq e^{C N},
\end{align} then
$$\lim_{\delta\rightarrow 0}\limsup_{N}P_N\left[\sup_{s,t}||\int_s^t \Psi(u)du||_{-m}> \epsilon\right] = 0.$$
\end{prop}
Hence it remains to prove \reff{tightone:line3}. We follow a
similar argument presented in \cite{Y92}. By the Cauchy-Schwarz
inequality
$$N||\int_s^t\Psi_udu||_{-m}(t-s)^{-1/2}\leq \frac{a^2N}{t-s}||\int_s^t\Psi_udu||^2_{-m}+\frac{N}{a^2}$$
and hence to obtain \reff{tightone:line3} it suffices to prove
that for a constant $C=C(a)>0$, independent of $t$ and $s$,
\begin{align}\label{tightone:line4}
E_{EQ}\left[\exp\left\{N\frac{a^2}{(t-s)}||\int_s^t\Psi_udu||^2_{-m}\right\}\right]\leq
e^{CN}.
\end{align}
Define $Q$ to be a Gaussian measure defined on $\mc H_{\alpha}$
with zero mean and covariance $S=(1-\triangle_N)^{-\beta}$, where
$\triangle_N f = N^2 \left[f(x+1/N)-2f(x)+f(x-1/N)\right]$.   The
existence of such a measure is guaranteed for $\beta>1/2$ by the
theory of Gaussian measures on Hilbert spaces (see for example
\cite{V68}). Define the inner product \linebreak $<F,G>_{k} =
<F,(1-\triangle_N)^k G>$. We thus have that
\begin{eqnarray*}
\int \exp \{<F,G>\} dQ(G) &=& \int \exp
\{<F,(1-\triangle_N)^{-\alpha} G>_{\alpha}\}dQ(G)\\
&=& \exp\left\{\frac{1}{2}||F||_{-\alpha-\beta}\right\}.
\end{eqnarray*}
We choose $m = \alpha + \beta$ and changing the order of
integration we may write \reff{tightone:line4} as
\begin{align}\label{tightone:line5}
\int E_{EQ}\left[\exp\left\{a
N^{1/2}(t-s)^{-1/2}<\int_s^t\Psi(u)du,f>\right\} \right]dQ(f) \leq
e^{C(a) N}.
\end{align}
Suppose that we could show the following bound for some positive
integer $k$
\begin{align}\label{tightone:line6}
&E_{EQ}\left[\exp\left\{a
N^{1/2}(t-s)^{-1/2}<\int_s^t\Psi(u)du,f>\right\}\right] \leq e^{C
a^2\{||f||_k+1\}}.&
\end{align}
Plugging condition \reff{tightone:line6} into
\reff{tightone:line5} we would obtain that the left hand side of
\reff{tightone:line5} is bounded by
\begin{align}\label{tightone:line7}
\exp\{Ca^2\} \int \exp\{Ca^2 <f,f>_k\}dQ(f),
\end{align}
which, because $Q$ is a Gaussian measure, integrates to
\begin{eqnarray*}
\{\det\left(1- 2 C a^2 (1-\triangle_N)^{k-m}\right)\}^{-1/2},
\end{eqnarray*}
for $a$ sufficiently small and $m-k>\frac{1}{2}$.  Since this is
(quite) smaller than $e^{CN}$ for some positive constant $C$ we
have reduced the proof of \reff{tightone:line4} again to showing
\reff{tightone:line6}.  That is, we have proved the following
result.

\begin{prop}\label{usefultight2}
Suppose that there exists a constant $C>0$ and an $N_0$ such that
for all $N>N_0$
\begin{align}\label{tightone:line8}
E_{EQ}[\exp\{a N^{1/2}(t-s)^{-1/2}<f,\int_s^t\Psi(u)du>\}] \leq
e^{C a^2\{||f||^2_k+1\}}.
\end{align} then for $m-k>\frac{1}{2}$
\begin{align*}
E_{EQ}[\int_0^T \int_0^T \exp\left( N \frac{||\int_s^t
\Psi(u)du||_{-m}}{\sqrt{t-s}}\right)ds dt]\leq e^{C N}.
\end{align*}
\end{prop}
We thus need to show that \reff{tightone:line8} holds for
$\Psi^N_{i,t}$ when $i=2,3$.  To do this we will repeat the
arguments of the local equilibrium principle for both terms.
Again, to simplify notation, we will denote $\Psi^N_{i,t}$ simply
as $\Psi$. Both terms satisfy the same bounds.  As before, we
begin with the Feynman-Kac inequality,
\begin{align*}
\log E_{EQ}[e^{\gamma N^{1/2}<f,\int_s^t\Psi(u)du>}] \leq\ \
\delta N^{1/2} \sup_{h }\sup_{|y|\leq R}
\left\{E_{\nu_{K,y}}[\gamma||f||_{\infty}
\Psi(x)h]-N^{3/2}\frac{N}{K}D_{K,y}(\sqrt{f})\right\}.
\end{align*}
We apply the logarithmic Sobolev inequality given in Theorem
\ref{lem:logsob}, and the entropy inequality as in \reff{SS3} and
\reff{line:lep5} to bound the line above by
\begin{align*}
\delta  N^{1/2} \sup_{|y|\leq R}\left\{ \frac{N^{5/2}}{K^3} \log
E_{\nu_{K,y}}\left[\,\exp\{||f||_{\infty}\gamma K^3/N^{5/2}
\Psi(x) \}\right]\right\},
\end{align*}
where $\gamma=\frac{a}{\sqrt{t-s}}$. We next apply Lemma \ref{LLD}
and obtain that this is smaller than $C(l, N_0, a_0)
(||f||_{\infty}+||f||_{\infty}^2)a^2$ for all $a>a_0$, and $N >
N_0$.  Because we are on the unit circle $\bb T$ we have, using
the Sobolev inequality $\sup_{x\in\bb T}f(x)\leq \int_{\bb
T}|f'(x)|dx$, for some new constant $C(l)$
\begin{eqnarray*}
\log E_{EQ}[\exp\{a N^{1/2}(t-s)^{-1/2}<\triangle_N
f,\int_s^t\Psi(u)du>\}] &\leq& C(l) a^2\{1+||f||_3\}.
\end{eqnarray*}
Combining this with Remark \ref{usefultight} and Remark
\ref{usefultight2} we may conclude that for $i=2,3$
\begin{align}\label{tightone:line77}
E_{EQ}[\int_0^T \int_0^T \exp\left( N (t-s)^{-1/2}||\int_s^t
\Psi_i(u)du||_{-m}\right)ds dt]\leq e^{C N}.
\end{align}
This concludes the proof of condition (T2).

By Chebychev's inequality as well assumption (F2), to prove that
condition (T1) is satisfied it is enough to show that
\begin{eqnarray}
&\lim_{M\rightarrow\infty}P_N\left[\sup_t \sup_{||f||_{m}\leq
1}\{<\triangle_Nf, \int_0^t \tilde{Y}_u^Ndu
>\}>M\right]=\ \ 0&\label{tightlineaaaa}\\
&\hspace{-2cm}\mbox {and}\,\,\,\,\,\,\,E_{NEQ}\left[\sup_t
||M^N_t||_{-m}\right]<\ \ \infty&\label{tightlinebbbb}
\end{eqnarray}
We begin with the latter.  By Lemma \ref{TBboundtails2} it is
enough to prove that
$$E_{NEQ}\left[\sup_t <M_t^N,f>\right]<\infty.$$
This follows from Doob's inequality and the result of Section
\ref{section:QV}.

The case \reff{tightlineaaaa} will be a repeat of the argument for
$(T2)_A$.  We split the field $\tilde Y$ as in \reff{chooseK}, and
show that
\begin{eqnarray*}
&\lim_{M\rightarrow\infty}P_N\left[\sup_t \sup_{||f||_{m}\leq
1}\{<\triangle_Nf, \int_0^t\Psi_{i,u}^Ndu
>\}>M\right]=0.&
\end{eqnarray*}
for $i=1, \ldots, 6.$  For $i\neq2,3$ we first apply Chebychev's
inequality.  For $i=1$ we argue as in \reff{BGbound1}.  For $i=4$
we use Lemma \ref{lem:hnorm}.   Lastly, we use \reff{BGbound45}
and Lemma \ref{lem:hnorm} to handle $i=5,6$. For the two remaining
cases, it turns out that we have already done all of the work. We
again use the Garcia-Rodemich-Rumsey inequality, Proposition
\ref{GRR}. We need simply to change the roles of $\delta$ and
$\epsilon$.  Again, we first argue for a general field $\Psi$. Let
$A_0$ and $B_0$ be as in \reff{defA0} and \reff{defB0}.
From
\reff{ineqBA1} and \reff{ineqBA2} it follows that for any field
$\Psi$, if we can show that \begin{eqnarray*} P_{EQ}^N(B_0 >
\frac{T^2}{4}e^{N\epsilon /16 \sqrt{T}}) \leq e ^{-N C(\epsilon)},
\end{eqnarray*}
where $C(\epsilon)\rightarrow\infty$ as $\epsilon \rightarrow
\infty$, then
\begin{eqnarray}
&\lim_{\epsilon\rightarrow
\infty}\limsup_{N}P_N[\sup_{t}||\int_0^t \Psi(u)du||_{-m}>
\epsilon] = 0,&
\end{eqnarray}
By Chebychev's inequality, it remains to show that
$E_{EQ}[B_0]\leq e^{C N}$ for some positive and finite constant
$C$.  However, this is exactly the content of
\reff{tightone:line77}.  This concludes the proof of tightness.
Notice again that we have in fact proved the following.

\begin{prop}\label{usefultight3}
Suppose that $E_{EQ}[B_0]\leq e^{C N}$, with $B_0$ defined in
\reff{defB0}, for some positive and finite constant $C$, then
$$\lim_{\epsilon\rightarrow \infty}\limsup_{N}P_N\left[\sup_{t}||\int_0^t \Psi(u)du||_{-m}> \epsilon\right] = 0.$$
\end{prop}


\section{Colour Density Fluctuations \label{sec:neqtwo}}


In this section we prove Theorem \ref{thm:main}.  For the sake of
brevity, we shall write out the proof for the case when $k=2$; the
generalization to any value of $k$ is immediate. As for the
colour-blind case, a key ingredient in the proof is the
logarithmic Sobolev inequality.  We use the additional condition
(E) in order to have the logarithmic Sobolev inequality, and the
remainder of the work is valid without this assumption.

The general outline of the proof is the same as the approach of
the previous section.  We show that the sequence $P_N$ is tight.
We then analyze the martingales under possible weak limits of
$P_N$ and use them to identify the form of the martingales under
the limiting measure.  We finish by showing that there is only one
measure which solves the martingale problem.

It is simpler to work with the following version of
\reff{resultintro}, for any smooth test functions $f_i$
\begin{align}\label{useSDEgen}
<dY_{i,t},f_i>&=&<Y_{i,t},\partial_{\rho^i}\varphi_i
f''>dt+<Y_{j,t},\partial_{\rho^j}\varphi_i f''>dt+<dW_t^i,
\sqrt{\varphi_i}f'_i>
\end{align}

We begin by identifying the martingales.  For $\bm f$ in $H_{4}^2$
define $<\bm Y_t^N, \bm f>$ to the pair
$\{<Y_{t,1}^N,f_1>,<Y_{t,2}^N,f_2>\}$. Under $P_N$ we have that
\begin{eqnarray}\label{id:mart2N}
&<\bm M_t^N,\bm f>=<\bm Y_t^N,\bm f>-< \bm Y_0^N,\bm f>-\int_0^t
<\tilde{\bm{Y}}_s^N,\triangle \bm f>ds&
\end{eqnarray}
is a martingale, where $\tilde{\bm Y}$ denotes the coupled
fluctuation field created by the pair of fields $\sqrt{N}~(c_i(\bm
\eta(\cdot))~-~\varphi_i(\bm \rho(\cdot))).$

The outline of this section is as follows.  In Section
\ref{sec:BG2} we identify the drift martingale under the limit of
any converging subsequence of the measures $P_N$.  To do this we
prove the Boltzmann-Gibbs principle.  We next handle the quadratic
variation and show that there is a unique solution to the
martingale problem. We finish the section by proving that the
sequence of measures $P_N$ is tight. As in the previous section
this completes the proof of Theorem \ref{thm:main}.

\subsection{Identifying the Drift. \label{sec:BG2}}

Arguing as in the previous section in \reff{id:martN} and
\reff{id:mart}, to show that the limiting martingales
\reff{id:mart2N} are the same as the drift martingale of
\reff{useSDEgen}, it is enough to prove the following.

\begin{thm}[Boltzmann-Gibbs]\label{BG2}
For any $f$ in  $C^1(\bb T)$ and $i=1,2$
\begin{eqnarray*}
\limsup_{N\rightarrow\infty}E_{NEQ}\left[\left|\int_0^t
<f,\sqrt{N}(c_i(\bm{\eta}^N_s)-\varphi_i(\bm{\rho}_s)-\nabla\varphi_i(\bm{\rho}_s)\cdot(\bm\eta^N_s-\bm\rho_s))>ds\right|\right]\rightarrow0.
\end{eqnarray*}
\end{thm}

The first step of the proof is to replace $\bm \rho_t(x)$, the
weak solution of $\partial_t\bm\rho=\frac{1}{2}\nabla \bm D_2(\bm
\rho)\nabla \bm \rho$ with the solution to the discretized version
\begin{align} \label{discreteHSLtwo}
\partial_t\bm\rho=\frac{1}{2}\nabla_N \bm D_2(\bm \rho)\nabla_N \bm\rho
\end{align}
with initial conditions $\bm \rho_t(x)|_{t=0}=\bm \rho_0(x)$ for
$x$ in $\bb T^N\times \bb T^N$.  The difference between the two
solutions is of order $\frac{1}{N}$ \cite{RM} and hence does not
affect our calculations. Because of this fact and in order to
simplify notation, we continue to denote the solution of
\reff{discreteHSLtwo} as $\bm\rho$ in the remainder of this
section.

We proceed with a proof of the two main ingredients:  the local
equilibrium principle and ``equivalence of solutions". Notice the
similarity to the proofs of the previous section. This happens
because the jump rates as well as their expectations for the
individual types are of the form
$\varphi_i=\rho^i\frac{\varphi(\rho)}{\rho}$, exhibiting linear
behaviour only in the density $\rho^i$.

\subsubsection{Local Equilibrium Principle.}
For a function $g:\bb{N}_+\times\bb{N}_+\mapsto\bb{R}$ define
\begin{eqnarray*}
V(x)=V_{(g,K,R)}(x)&=&AV_{|y-x|\leq
K/N}\left\{g(\bm{\eta}^N(y))-E_{\mu_{\bm{\rho}=\bm{m}^K(x)}}[g(\bm{\eta}^N(y))]\right\}\bb{I}[m^K(x)\leq R]\\
&=&\left\{\bar{g}_K(\bm{\eta}^N(x))-\hat{g}_K(\bm{m}^K(x))\right\}\bb{I}[m^K(x)\leq
R],
\end{eqnarray*}
with $m^K(x)$ defined as in the previous section.
\begin{lem}\label{LEP2}
Suppose that $g$ is a function satisfying $|g(k_1,k_2)|\leq
C(1+k_i)$, for all $(k_1,k_2)$ and for some $i$, then
\begin{align*}
\lim_{l\rightarrow0}\lim_{N\rightarrow\infty}\sup_{||J||_{\infty}<1}N^{1/2}E_{NEQ}[|\int_0^t<J,V_s>ds|]\leq0.
\end{align*}
\end{lem}
\noindent Without loss of generality, we will consider only the
case $g(k_1,k_2)\leq C (1+k_1)$ in the proof.  We begin with a lemma.  
As in the previous section we state the next result for
configurations $\eta$ prior to the diffusive re-scaling.  For $g$
as in Lemma \ref{LEP2}, define $\tilde g = g - E_{ \mu_{\rho^i =
m^K_i(x), i=1,2}(\cdot|\eta_2)}[g],$ where $m^K_i(x)=
AV_{|x-y|\leq K} \,\eta_i(x)$.

\begin{lem}[Inhomogeneous Local Large Deviations]\label{LLD2}
There exists a positive constant $C$ independent of the
configuration $\eta_2$ such that
\begin{eqnarray*}
\log E_{\nu_{K,\bm y}(\cdot|\eta_2)}\left[\exp\left\{\gamma AV_{|x-y|\leq K} \tilde{g}(\bm{\eta}(x))\right\}\right]\leq C\{\frac{\gamma^2}{K}+\frac{\gamma+\gamma^2}{K^2}\}.\notag\\
\end{eqnarray*}
\end{lem}
\begin{proof}
We first note that $\eta_1$ conditioned on $\eta_2$ satisfies the
conditions of a nonhomogeneous zero-range model as was shown in
Section \ref{relateZR}, and hence we have the necessary local
limit result from Theorem \ref{normalLLT} with bounds which are
independent of $\eta_2$ (and hence the site $x$). We then repeat
the same argument as in the proof of Lemma \ref{LLD}.
\end{proof}

\begin{proof}[Proof of Lemma \ref{LEP2}.]
Let $X(s)$ denote $<J,V_s>$, as before. 
Most of the proof is the same here as in the colour-blind case,
and hence we omit many of the details.  By the entropy inequality
and Feynman-Kac formula, we have
\begin{align}
E_{NEQ}[N^{1/2}|\int_0^tX(s)ds|]
 \leq \frac{2}{\beta N} t \sup_s
 \Gamma^N(s)+\frac{C}{\beta},\label{lineSSS22}
 \end{align}
where $\Gamma^N(s)$ is the largest eigenvalue of $\pm \beta
N^{3/2}X(s) + L_N$.  The biggest difference in this setting is
that we need to introduce an additional conditioning.  We reduce
the span of the Dirichlet form as in \reff{SS2} and condition on
the average density \emph{and} the particles of colour two to
obtain that $\Gamma^N(s)$ is bounded above by
\begin{align*}
N^{3/2} \sup_{f} \sup_{|\alpha|\leq\beta||J||_{\infty}}\sup_{|\bm
y|\leq \bm
R}\sup_{\eta_2}\{E_{\tilde{\nu}_{\Lambda_K(x),y_1}}[\alpha \cdot V
\cdot f ]- \frac{N^{3/2}}{2K+1} D_{\nu_{K,\bm
y}(\cdot|\eta_2)}^{K,x,1}(\sqrt{f})\}
\end{align*}
where $D^{x,1}_{\nu_{K,y_1,y_2}}(f)$ defined in
\reff{conddirform}. Using the logarithmic Sobolev inequality
(Theorem \ref{thm:logsobnh}), the entropy inequality and Lemma
\reff{LLD2} as in the colour-blind proofs we conclude that
\begin{align*}
E_{NEQ}[N^{1/2}|\int_0^tX(s)ds|] \leq
C\{\beta^2l^2+\frac{l^{-2}+\beta^3l^4}{N^{1/2}}\}+\frac{C}{\beta},
\end{align*}
from which the result follows.
\end{proof}

\subsubsection{Equivalence of Solutions.}

Define $\bar{Y}^N_{i,t}$ to denote the centered version of the
function ${Y^N_{i,t}}$ as in \reff{def:Ybar}.  Because of the
symmetry of the process, we only consider
$||\bar{Y}^N_{1,t}||_{*}^2$. The same formula holds as in the
colour-blind case:
\begin{eqnarray}
&&\hspace{-1cm}E_{NEQ}[N^{-1/2}||\bar{Y}^N_{1,t}||_{*}^2+2<\bar{\eta}^N_1-\bar{\rho}_1,\bar{Y}^N_{1,t}>_{*}]\notag\\
&&\hspace{2cm} - \ \ E_{NEQ}[N^{-1/2}||\bar{Y}^N_{1,0}||_{*}^2+2<\bar{\eta}^N_1-\bar{\rho}_1,\bar{Y}^N_{1,0}>_{*}]\notag\\
   &=&  E_{NEQ}[\int_0^t \{-<c_1(\bm\eta^N)-\varphi_1(\bm{\rho}),\bar{Y}^N_{1,s}> +
  N^{-1/2}(1+1/N)\sum_{x\in\bb{T}^N}c_1(\bm{\eta}^N_s(x))\}ds]\notag\\
  &=& E_{NEQ}[\int_0^t N^{1/2}AV_{x\in\bb{T}^N} H_{1,s}(x)ds]+N^{-1/2}E[\int_0^t AV_{x\in\bb{T}^N}c_1(\bm{\eta}^N_s(x))ds],\label{line:2hnorm3}
\end{eqnarray}
where
$H_1(x)=c_1(\bm{\eta}^N(x))-(c_1(\bm{\eta}^N(x))-\varphi_1(\bm{\rho}(x)))(\eta_1^N(x)-\rho^1(x))$.

We next replace $H_1(x)$ with
$$\bar{H}^K_1(x)=AV_{|x-y|\leq
K/N}[c_1(\bm{\eta}^N(y))-\{c_1(\bm{\eta}^N(y))-\varphi_1(\bm{\rho}(x))\}\{\eta_1^N(y)-\rho^1(x)\}].$$
By integration by parts and \reff{SSS} we conclude that the error
produced by this replacement is of the order of $O(l)$. Line
\reff{line:2hnorm3} is thus equal to
\begin{align}
E_{NEQ}[\int_0^t N^{1/2}AV_{x\in \bb T^N}\bar{H}^K_{1,s}(x)ds] +
  O(l)+O(N^{-1/2}).\label{line:2hnorm4}
\end{align}
Next define
$F_1(\bm{\eta}^N(y),\bm{\rho}(x))=c_1(\bm{\eta}^N(y))-\{c_1(\bm{\eta}^N(y))-\varphi_1(\bm{\rho}(x))\}\{\eta_1^N(y)-\rho^1(x)\}$,
and split $\bar{H}^K(x)$ into three pieces:
\begin{eqnarray*}
\bar{H}^K_i&=&AV_{|x-y|\leq K/N}[F_1(\bm{\eta}^N(y),\bm{\rho}(x))-\hat{F}_1(\bm{m}^K(x),\bm{\rho}(x))]\bb{I}(m^K\leq R)\bb{I}(s^K_1(x)\leq R_1)\\
&&+AV_{|x-y|\leq K/N}[F_1(\bm{\eta}^N(y),\bm{\rho}(x))]\bb{I}(m^K > R\cup s^K_1(x)> R_1)\\
&&+AV_{|x-y|\leq K/N}[\hat{F}_1(\bm{m}^K(x),\bm{\rho}(x))]\bb{I}(m^K\leq R)\bb{I}(s^K_1(x)\leq R_1)\\
&=&AV_{|x-y|\leq K/N}[F_1^{(1)}(x)+F^{(2)}_1(x)+F^{(3)}_1(x)],
\end{eqnarray*}
where
\begin{eqnarray*}
\hat{F}_1(\bm{m}^K(x),\bm{\rho}(x))&= &E_{\mu_{\bm{\rho}=\bm{m}^K(x)}}[F_1(\bm{\eta}^N(y),\bm{\rho}(x))]\\
&=&-(\varphi_1(\bm{m}^K(x))-\varphi_1(\bm{\rho}(x)))(m^K_1(x)-\rho^1(x))],
\end{eqnarray*}
and $s^K_1(x)=AV_{|x-y|\leq K/N}\,\{\eta_1^N(x)\}^2$.

By the local equilibrium principle the first term vanishes in the
limit.  Even though the term $F$ is not linear in $\eta$ as
required in the local equilibrium principle, the addition of the
condition $\bb{I}(s^K_1(x)\leq R_1)$ fixes the problem originally
pointed out in Remark \ref{rem:lepfem}. By an identical argument
to the one presented in the previous section, for fixed $R$ and
$R_1$ large enough there exists a constant $\tilde{A}_2>0$ such
that the following bound holds:
\begin{align}
AV_{|x-y|\leq K/N}F^{(2)}_1(x)\leq
-\tilde{A}_2(m^K_1(x)-\rho^1(x))^2\bb{I}(m^K(x)> R\cup
s^K_1(x)>R_1)\label{line:2hnorm6}.
\end{align}
We handle $F^{(3)}_1$ similarly.  In this case we make use of the
fact that the function $\varphi_1(\rho)$ is strictly increasing
(Proposition \ref{ch1:prop:lip}).  We thus have that there exists
an $\tilde{A}_3>0$ such that:
\begin{align}
AV_{|x-y|\leq K/N}F^{(3)}(x)
\leq -\tilde{A}_3(m^K_1(x)-\rho^1(x))^2\bb{I}(m^K(x)\leq
R)\bb{I}(s^K_1(x)\leq R_1)\label{line:2hnorm7}.
\end{align}
Combining the arguments from \reff{line:2hnorm4} through
\reff{line:2hnorm7}, we find that
\begin{align}
&E_{NEQ}[N^{-1/2}||\bar{Y}^N_{t,1}||_{*}^2+2<\bar{\eta}_1-\bar{\rho}_1,\bar{Y}^N_{t,1}>_{*}]- E_{NEQ}[N^{-1/2}||\bar{Y}^N_{0,1}||_{*}^2+2<\bar{\eta}_1-\bar{\rho}_1,\bar{Y}^N_{0,1}>_{*}]\notag\\
  &\leq-\tilde{A}_2\wedge \tilde{A}_3E_{NEQ}[\int_0^t N^{-1/2}\sum_{x\in \bb T^N}(m_1^K(x)-\rho^1(x))^2 ds] +
  O(l^{-1})+o(1)
  \label{line:2hnorm8}
\end{align}
Using the bound \reff{eos:line} together with assumptions (F1) on
\reff{line:2hnorm8} as in the colour-blind proof, we conclude
\begin{thm}\label{normtrickresult:ktype}
\begin{align*}
\lim_{l\rightarrow0}\lim_{N\rightarrow\infty}E_{NEQ}[\int_0^t
N^{-1/2}\sum_{x\in \bb T^N}(m^K_i(x)-\rho^i(x))^2ds]=0,
\,\,\,i=1,2.
\end{align*}
\end{thm}

We are now ready to prove Theorem \ref{BG2}.  Define the field
$\bm{\Phi}^N(t)=\{\Phi^N_1(t),\Phi_2^N(t)\},$ with $\Phi^N_1(t)$
given by
\begin{eqnarray*}
\sqrt{N}\left\{c_1(\bm{\eta}^N_t)-\varphi_1(\bm{\rho}_t)
-\partial_1\varphi_1(\bm{\rho}_t)(\eta_{1,t}^N-\rho_{1,t})
-\partial_2\varphi_1(\bm{\rho}_t)(\eta_{2,t}^N-\rho_{2,t})\right\},
\end{eqnarray*}
and $\Phi^N_2(t)$ defined similarly.

As before, we split $\bm\Phi^N$ into five separate parts.
\begin{eqnarray*}
 \bm{\Phi}^N(t)=\bm{\Phi}^N_1(t)+\bm{\Phi}^N_2(t)+\bm{\Phi}^N_3(t)+\bm{\Phi}^N_4(t)+\bm{\Phi}^N_5(t),
\end{eqnarray*}
where the first entry in each is given by
\begin{eqnarray*}
N^{-1/2}\Phi^N_{1,1}&=&c_1(\bm{\eta}^N(x))-\bar{c}_1^K(\bm{\eta}^N(x))\\
        &
        &\hspace{1cm}-\partial_1\varphi_1(\bm{\rho}(x))\{\eta^N_1(x)-m^K_1(x)\}-\partial_2\varphi_1(\bm{\rho}(x)\{\eta^N_2(x)-m^K_2(x)\}\\
N^{-1/2}\Phi^N_{1,2}&=&\{\bar{c}_1^K(\bm{\eta}^N(x))-\varphi_1(\bm{m}^K(x))\}\bb{I}[m^K(x)\leq R]\\
N^{-1/2}\Phi^N_{1,3}&=&\left[\varphi_1(\bm{m}^K(x))-\varphi_1(\bm{\rho}(x))\right.\\
            &&\hspace{1cm}-\partial_1\varphi_1(\bm{\rho}(x))\{m^K_1(x)-\rho^1(x)\}\\
            &&\hspace{1cm}\left.-\partial_2\varphi_1(\bm{\rho}(x))\{m^K_2(x)-\rho^2(x)\}\right]\bb{I}[m^K(x)\leq R]\\
N^{-1/2}\Phi^N_{1,4}&=&\left[-\varphi_1(\bm{\rho}(x))\right.\\
            &&\hspace{1cm}-\partial_1\varphi_1(\bm{\rho}(x))\{m^K_1(x)-\rho^1(x))\}\\
            &&\hspace{1cm}\left.-\partial_2\varphi_1(\bm{\rho}(x))\{m^K_2(x)-\rho^2(x)\}\right]\bb{I}[m^K(x)> R]\\
N^{-1/2}\Phi^N_{1,5}&=&\bar{c}_1^K(x)\bb{I}[m^K(x)> R],
\end{eqnarray*}
with $\bar{c}_1^K(x)=AV_{|x-y|\leq\frac{K}{N}}c_1(\eta(y))$. By
symmetry, we only work with the $\Phi^N_1(t)$. By the local
equilibrium principle and Theorem \ref{normtrickresult:ktype} we
already have for $i=2$ and 3 that
\begin{align}\label{toprove}
E_{NEQ}[|\int_0^T<\Phi_{1,i}^N,f>ds|]\rightarrow 0,
\end{align}
for any continuous function $f$.    The result for $i=5$ follows
from the same argument for the colour blind case, because
$\bar{c}_1^K(x)\leq\bar{c}^K(x)$.  Next, we consider $i=4$. We may
pick $ R$ sufficiently large so that
\begin{eqnarray*}
|\Phi^N_4(t)|\leq C N^{-1/2}\sum_{i=1,2}^N(m^K_i(x)-\rho^i(x))^2
\end{eqnarray*}
and hence \reff{toprove} holds also for $i=4$.  It remains to
consider $i=1$.  An argument using summation by parts like the one
presented in the previous section proves the result. Notice that
it is in this part of the argument that we require the most
smoothness of the function $f$.  It is enough to consider test
functions $f$ which have a bounded first derivative. This
concludes the proof of the Boltzmann-Gibbs principle.

\subsection{Quadratic Variation.}

We begin with a lemma.  Let $K<<N$, fix $R$ large, and define
$$X_{t,i}^K(x)=AV_{|x-y|\leq \frac{K}{N}}\{ c_i(\bm \eta^N_t(x))-\varphi_i(\bm\rho_t(x))\}\bb I [m^K(x)\leq R].$$

\begin{lem}\label{expbound}
Assume that $|t-s|<T$. For any choice of $l=\frac{K}{\sqrt{N}}$
there exists a constant $C=C(l,R)$ so that the following bound
holds uniformly for $||J||_{\infty}\leq 1$:
$$\log E_{EQ}[\exp\{\gamma N <J,\int_s^tX^K_{u,i}du>\}]\leq C (\gamma^2 +\gamma)T.$$
The same bound holds if we replace $J$ by a function dependent on
time with uniform bound $||J_t||_{\infty}\leq 1$.
\end{lem}

\begin{proof}
By the Feynman-Kac inequality, and conditioning as before, we have
that the left hand side is bounded above by
\begin{align*}
& \,\,T N \sup_{f\geq 0: E_{\mu_{\rho^1,\rho^2}}[f]=1}
\{E_{\mu_{\bm\rho}}[\gamma
X_if]]-ND_{N,\bm\rho}(\sqrt{f})\}\\
\leq&\,\,T N \sup_{f }\sup_x \sup_{|\bm y|\leq \bm R}
\sup_{\eta_2} \{E_{\nu_{K,\bm y}(\cdot|\eta_2)}[\gamma
X_i(x)f]]-\frac{N^2}{K}D^{x,1}_{\nu_{K,\bm
y}(\cdot|\eta_2)}(\sqrt{f})\}
\end{align*}
We next apply the logarithmic Sobolev inequality, Theorem
\ref{thm:logsobnh}, followed by the entropy inequality to obtain
that this is less than
\begin{align*}
T \frac{N^3}{K^3} \sup_{|\bm y|\leq \bm R}\sup_{\eta_2}\{\log
E_{\nu_{K,\bm y}(\cdot|\eta_2)}[\exp\{\gamma K^3/N^2 X_i\}]\}.
\end{align*}
Using Lemma \ref{LLD2} in the above along with $l\sqrt{N}=K$
proves the result.
\end{proof}

We now turn to the calculation of the quadratic variation.  For
fixed $\alpha_1$ and $\alpha_2$ we have that $\alpha_1
<M_{t,1}^N,f_1>+\alpha_2 <M_{t,2}^N,f_2>$ has quadratic variation
under $P_N$ given by:
\begin{eqnarray}\label{martintwo2}
\hspace{-.5cm}\sum_i\frac{\alpha_i^2}{2}\int_0^t N^{-1}\sum_{x\in
       \bb{T}^N}[c_i(\eta^N_{s}(x))+c_i(\eta^N_{s}(x+\frac{1}{N})][\nabla_Nf(x)]^2ds.
\end{eqnarray}
To show that in the limit these martingales match with
\reff{useSDEgen} it is enough to prove the following proposition
holds for $i=1,2$.
\begin{prop}\label{QVtwo}
For a smooth function $h$ in $C^1(\bb T)$.
\begin{align*}
E_{NEQ}\left[\left|\int_0^t<h,c_i(\eta^N_s)-\varphi_i(\rho_s)>ds\right|\right]\rightarrow
0 .
\end{align*}
\end{prop}

\begin{proof}
Let $X_{s,i}$ denote the field $c_i(\eta^N_s)-\varphi_i(\rho_s)$
for $i=1,2$. Our first step is to replace $X_{s,i}$ with an
average over boxes of size $K$.  By summation by parts, the
replacement adds a term of size $\frac{K}{N}$.  Hence we need only
worry about the inequality using the field averaged over boxes of
size K : $X^K_{s,i}$.  Fix a large $R$ and consider
\begin{eqnarray*}
X^K_{s,i}(x)
        &=& AV_{|x-y|\leq K/N}\left[c_i(\eta^N_s(x))-\varphi_i(\rho_s(x))\right]\bb I[m^K(x)\leq R]\\
        &&+ AV_{|x-y|\leq K/N}\left[c_i(\eta^N_s(x))-\varphi_i(\rho_s(x))\right]\bb I[m^K(x)> R]\\
        &=& X^{K,(a)}_{s,i}(x)+X^{K,(b)}_{s,i}(x)
\end{eqnarray*}
By Chebychev's inequality the second field, $X^{K,(b)}_i(x)$, is
bounded by
\begin{eqnarray*}
E_{NEQ}[\int_0^t |<h,X^{K,(b)}_{i,s}>|ds]\leq \frac{C}{R}
\end{eqnarray*}
for some constant C.  By the entropy inequality, it follows from
Lemma \ref{expbound} and assumption (H1) that
\begin{eqnarray*}
E_{NEQ}[\int_0^t <h,X^{K,(a)}_{i,s}>ds] \leq  \frac{C \beta
}{N}+\frac{C}{\beta }.
\end{eqnarray*}
The result follows if we first let $N$ converge to $\infty$, then
do the same for $\beta$, then $R$.
\end{proof}

\subsection{Uniqueness of the Martingale Problem.}

In this section we show that the martingale problem described by
the generalized stochastic differential equation
\reff{resultintro} has a unique solution.  We state and prove the
theorem only for 2 colours, however, the statement as well as the
proof are valid for any number of colours.

Define for $\bm f$ in $\mc H_{m}^2$
\begin{eqnarray}\label{line:marttwo}
\hspace{-1cm}<\bm M_t,\bm f>&=&<\bm Y_t,\bm f>-<\bm Y_0,\bm
f>-\int_0^t <\frac{1}{2}\bm D_2(\bm \rho)^* \bm{Y}_s,\triangle \bm
f>ds.
\end{eqnarray}
We also define the operators $\mc S_{t,i}=\sqrt{\varphi_i(\bm
\rho_t )}\nabla$, $\mc A_t^i=\sqrt{S'(\rho_t)\rho_{t,i}
}\triangle$ and $\mathcal{P}^S_{t,u}$ the semigroup associated to
the operator $\frac{1}{2}S(\rho_t)\triangle.$

\begin{thm}\label{uniquetwo}
Fix a positive integer $m\geq2$.  Let $P$ be a probability measure
on the space $\{(C[0,T],H_{-m}^2),\mathcal{ F}\}$, where
$\mathcal{F}=\cup_{t\geq0}\mathcal{F}_t$, and $\mathcal{F}_t$ is
the canonical filtration of the process $<\bm Y_t,\bm f>$. Assume
that for any $\bm f \in C^{\infty}(\mathbb{T})\times
C^{\infty}(\mathbb{T}) $ and $\bm \alpha\in \bb R^2$ the
quantities $<\bm M_t,\bm f>$ defined in \reff{line:marttwo}, and
\begin{eqnarray}\label{line2:marttwo}
&\{\bm \alpha \cdot <\bm M_t,\bm f>\}^2-\sum_{i=1,2} \alpha_i^2
\int_0^t||\mathcal{S}_{s,i}f_i||^2_2ds
\end{eqnarray}
are $L^1(P)$ continuous martingales with respect to
$\mathcal{F}_t$.  Then, for any colour $i$, for all $0\leq s\leq
t$, $f$ in $C^{\infty}(\mathbb{T})$, and subsets $A$ of
$\mathbb{R}$,
\begin{eqnarray}
&&\hspace{-2cm}P[\left\{ <Y_{t,i},f>-\int_s^t <Y_u,\mc A_u^i \mc P^S(t,u)f>du \right \}\in A|\mathcal{F}_s]\notag\\
&= &\int_A \frac{1}{\sqrt{2\pi \int_s^t
||\mathcal{S}_{u,i}\mathcal{P}^S_{t,u}f||^2_2du}}\exp\left\{-\frac{(y-<Y_{i,s},\mathcal{P}^S_{t,s}f>)^2}{2\int_s^t
||\mathcal{S}_{u,i}\mathcal{P}^S_{t,u}f||^2_2du}\right\}dy,\label{uniqueinvert}
\end{eqnarray}
$P$ almost surely.
\end{thm}

From this result we obtain uniqueness of the limiting measure $P$.
From the above theorem we not only have the distribution of
$$<Y_{t,i},f>-\int_s^t <Y_u,\mc A_u^i \mc P^S_{t,u}f>$$
for any colour $i$, but we also know that these quantities are
independent for different $i$ and different test functions $f$
because of \reff{line2:marttwo}.  If we choose the same $f$ and
sum over the colours we obtain
$$<Y_t,f>-\int_s^t <Y_u, \mc{A}_u \mc P^S_{t,u}f>,$$
where $\mc A_u$ is the sum over $\mc A_u^i$.  This fact plus
Theorem \ref{uniqueone} which gives the full distribution of the
process $<Y_t,f>$ is enough to obtain uniqueness of the pair
$\{<Y_t^1,f>,<Y_t^2,h>\}$. A standard Markov argument gives the
uniqueness of the finite dimensional distributions, which implies
the uniqueness of the measure $P$.

The proof of the above theorem is similar to the proof of the
uniqueness result for the colour-blind density fluctuation field
as it appears in \cite{K-L} or \cite{Ho-S}.

\begin{proof}
First of all notice that the martingale $<M_{t,i},f>$ in
\reff{line:marttwo} may also be written as
\begin{eqnarray*}
<M_{s,i},f>+<Z_{t,s}^i,f> - \int_s^t <Y_u,\mc A_u^i f>du,
\end{eqnarray*}
where $<Z_{t,s}^i,f>$ is used to denote
$$<Y_{t,i},f>-<Y_{s,i},f>-\int_s^t <Y_{u,i},\frac{1}{2} S(\rho_u)\triangle
f>du.$$

 It\^{o}'s formula together with the fact that \reff{line:marttwo} and \reff{line2:marttwo}  are both
 martingales imply that for each fixed $s\geq 0$ and $f \in
C^\infty(\bb T)$, $X_{t,s}(f)$ defined by
\begin{eqnarray*}
X_{t,s}(f)=\exp\left\{i\left(<Z_{t,s}^i,f> - \int_s^t <Y_u,\mc
A_u^i
f>du\right)+\frac{1}{2}\int_s^t||\mathcal{S}_{u,i}f||^2_2du\right\}
\end{eqnarray*}
is a martingale.  Next take two fixed times $t_1<t_2\leq T$ and
let $s_{n,j}$ denote a partition of the time interval defined by
$s_{n,j}=t_1+\frac{j}{n}(t_2-t_1).$  Consider the quantity defined
as
$$\Pi_{j=0}^{n-1}X_{s_{n,j+1},s_{n,j}}(\mc  P^S(T,s_{n,j})f).$$
By continuity of $<Y_s^i,\mc P^S(T,t)f>$ in $\{s,t\}$ we may let
$n \rightarrow \infty$ to obtain that the above quantity converges
a.s. and in $L^1(P)$ to $\frac{Z_{t_2}(f)}{Z_{t_2}(f)},$ where
$Z_t(f)$ is equal to
\begin{eqnarray*}
\exp\left\{i <Y_{t,i},\mc P^S_{T,t}f>-\int_0^t <Y_u,\mc A_u^i\mc
P^S_{T,u} f>du +\frac{1}{2}\int_0^t||\mathcal{S}_{u,i}\mc
P_{T,u}f||^2_2du\right\}.
\end{eqnarray*}
Because the convergence above takes place also in $L^1(P)$ the
martingale properties of $X_{t,s}(\mc P^S_{(T,s)}f)$ are passed
onto $Z_{t}(f)$. The fact that $Z_t(f)$ is a martingale proves the
theorem.
\end{proof}

\subsection{Tightness}

Let $P_N$ denote the probability measure on $D([0,T],
\mathcal{H}_{-m}^2)$ induced by the fluctuation field $\bm Y_t^N$.

For a function $\mb F=\{F^1,F^2\}$ in $D([0,T],
\mathcal{H}_{-m}^2)$, define the (uniform) modulus of continuity,
$\omega_\delta(\mb F)$, for a fixed $\delta>0$:
\begin{eqnarray*}
\omega_\delta(\bm F)=\sup_{|s-t|\leq\delta, 0\leq s,t \leq
T}\left\{||F^1_t-F^1_s||_{-m}+||F^2_t-F^2_s||_{-m}\right\}.
\end{eqnarray*}
As previously we shall simplify the supremum in the above notation
to $\sup_{s,t}$.  By Prohorov's theorem and  Arzela-Ascoli it
follows that to show the sequence $P_N$ is tight we need to show
that it satisfies two conditions:
\begin{enumerate}
\item[](T1)\hspace{0.5cm}$\lim_{B \rightarrow \infty}\limsup_{N \rightarrow \infty}
P_N [\sup_{0\leq t \leq T} ||\bm Y_t||_{-m}>B] =0$
\item[](T2)\hspace{0.5cm}$\lim_{\delta \rightarrow 0}\limsup_{N \rightarrow
\infty}P_N[\omega_\delta(\bm Y)> \epsilon]=0
\,\,\,\,\,\,\,\,\,\,\,\, \forall \epsilon
>0.$
\end{enumerate}

\begin{prop}\label{prop:tighttwo}
The sequence of measures $P_N$ is tight in $D([0,T],
\mathcal{H}_{-4}^2)$. Moreover, all limit points are concentrated
on continuous paths.
\end{prop}
Using \reff{id:mart2N} we have that
\begin{align*}
<Y_{t,1}^{N},f> \, = \,
<Y_{0,1}^N,f>+\int_0^t<\tilde{Y}_{s,1}^N,\triangle_Nf>ds+<M_{t,1}^N,f>,
\end{align*}
where $M_{t,1}^N$ is a martingale. By assumption (F2) there is
nothing to prove for the initial field. Also, $M_{t,1}^N$ has
quadratic variation given in \reff{martintwo2} with $\alpha_1=1$
and $\alpha_2=0$. This formula is very similar to the one in the
single colour case, and hence, the proof will follow as before
from the work we have already done in the two-colour
Boltzmann-Gibbs principle.  We may also handle each of $Y_{t,i}^N$
separately.

To show that $\bm Y^N_{t}$ satisfies condition (T2) it is enough
to show
\begin{enumerate}
\item[] $(T2)_A$\hspace{0.5cm}$\lim_{\delta\rightarrow 0}\limsup_{N\rightarrow\infty}P_N[\sup_{s,t}\sup_{||f||_{m}\leq 1}\{<\triangle_Nf, \int_s^t \tilde{Y}_u^{N,i}du
>\}> \frac{\epsilon}{4}]$
\item[] $(T2)_B$\hspace{0.5cm}$\lim_{\delta\rightarrow 0}\limsup_{N\rightarrow\infty}P_N[\sup_{s,t} ||M_{t,i}^N-M_{s,i}^N||_{-m} >
\frac{\epsilon}{4}]$,
\end{enumerate}
for each of $i=1,2$.

We begin with $(T2)_B$, and fix $i=1$.  Clearly, the other case is
the same.  Because of this special ``decoupling" that occurs, the
proof is now almost identical to the colour-blind case, and we
omit many of the details. The next three lemmas follow from the
quadratic variation result in Proposition \ref{QVtwo}, as before;
see the proof of Lemmas \ref{TBboundtails}, \ref{TBboundtails2}
and \ref{martboundsone}.

\begin{lem}[Bounding the tails 1.]\label{TBboundtails21}
There exists a finite constant $C=C(\varphi)$ such that for every
eigenfunction $e_z$ of $\mc H_m$,
$$\limsup_{N\rightarrow\infty}E_N[\sup_{0\leq t\leq T}|<M_{t,1}^N,e_z>|^2]\leq C T \{1+<\triangle_N e_z,\triangle_N e_z>\}$$
\end{lem}

\begin{lem}[Bounding the tails 2.]\label{TBboundtails22}
For $m\geq 4$,
$$\lim_{K\rightarrow\infty}\limsup_{N\rightarrow\infty}E_N[\sup_{0\leq t\leq T}\sum_{z\geq K}\{<M_{t,1}^N,e_z>\}^2\gamma^{-m}_z]=0.$$
\end{lem}

\begin{lem} For any $f$ in $C^3(\bb T)$
$$\lim_{\delta\rightarrow
0}\limsup_{N\rightarrow\infty}P_N[\sup_{s,t}\{<M_{t,1}^N-M_{s,1}^N,f>\}^2>\epsilon']=0.$$
\end{lem}

These three lemma together prove that condition $(T2)_B$ is
satisfied.  We next turn to condition $(T2)_A$: yet again we split
this field up in the sum of several terms:
\begin{eqnarray}
\Psi^{N,1}_{1,t}&=&\sqrt{N}\{c_1(\bm \eta^N_t(x))-\bar{c}_{K,1}(x)\}\notag\\
\Psi^{N,1}_{2,t}&=&\sqrt{N}\{\bar{c}_{K,1}(\bm \eta^N_t(x))-\varphi_1(\bm m^K_t(x))\}\bb{I}[m^K_t(x)\leq R]\notag\\
\Psi^{N,1}_{3,t}&=&\sqrt{N}\partial_1\varphi_1(\bm \rho_t(x))\{m^K_{1,t}(x)-\rho_{1,t}(x)\}\bb{I}[m^K_t(x)\leq R]\notag\\
\Psi^{N,1}_{4,t}&=&\sqrt{N}\partial_2\varphi_1(\bm \rho_t(x))\{m^K_{2,t}(x)-\rho_{2,t}(x)\}\bb{I}[m^K_t(x)\leq R]\notag\\
\Psi^{N,1}_{5,t}&=&\sqrt{N}\left[\varphi_1(\bm m^K_t(x))-\varphi_1(\bm \rho_t(x))\right.\label{chooseKtwo}\\
                && \hspace{1cm}\left.-\sum_{i=1}^2\partial_i\varphi_1(\bm \rho_t(x))\{m^K_{i,t}(x)-\rho_{i,t}(x)\}\right]\bb{I}[m^K_t(x)\leq R]\notag\\
\Psi^{N,1}_{6,t}&=&-\sqrt{N}\varphi_1(\bm \rho_t(x))\bb{I}[m^K_t(x)>\tau]\notag\\
\Psi^{N,1}_{7,t}&=&\sqrt{N}\bar{c}_{K,1}(x)\bb{I}[m^K_t(x)>\tau].\notag
\end{eqnarray}
We need to prove that
\begin{eqnarray}\label{tighttoprove}
\lim_{\delta\rightarrow
0}\limsup_{N\rightarrow\infty}P_N\left[\sup_{s,t}\sup_{||f||_{m}\leq
1}\left|\int_s^t \Psi_{i,u}^{N,1}(\triangle_Nf)du \right|>
\frac{\epsilon}{28}\right],
\end{eqnarray}
for $i=1, \ldots 7$.  As before, the quantities $i=2,3,4$ are more
difficult.

For $i=1,5,6,7$ we have that \reff{tighttoprove} holds because we
have already shown in the previous section
\begin{eqnarray*}
\lim_{l\rightarrow
0}\lim_{N\rightarrow\infty}E_{N}\left[\sup_{||f||_{m}\leq
1}\int_0^T\left| <\Psi_{i,u}^{N,1},\triangle_Nf>\right|du
\right]=0.
\end{eqnarray*}
That is, the result follows by summation by parts for $i=1$, and
by Theorem \ref{normtrickresult:ktype}  and $|\Psi_{i,u}^N(x)|\leq
C \sum_{i=1,2}(m_i^K(x)-\rho^i(x))^2$ for the remaining cases.

It remains to show \reff{tighttoprove} for $i=2,3,4$.  By
Propositions \ref{usefultight} and \ref{usefultight2} this follows
if we have the following bound for $i=2,3,4$
\begin{eqnarray*}
E_{EQ}]\left[\exp\{a N^{1/2}(t-s)^{-1/2}<\triangle_N
f,\int_s^t\Psi_{i,u}^Ndu>\}\right] \leq e^{C a^2\{||f||^2_k+1\}},
\end{eqnarray*}
for a positive constant $C$ and all large enough $N$.  As before,
to obtain convergence in $m\geq 4$ we need $k=3$.  Indeed this
follows from Lemma \ref{expbound} and the Sobolev inequality for
$i=2$. It is straightforward to check that the exact same bounds
apply for $i=3,4$.  The constant $C$ depends on $l$ and is valid
for all $a>a_0$, that is, $C=C(l,a_0)$.  This completes the proof
of $(T2)_A$.

To complete the proof of tightness we hence need the bound (T1).
Because of Lemma \ref{TBboundtails21} and \ref{TBboundtails22} we
have that for $i=1,2$
$$\lim_{M \rightarrow \infty}\limsup_{N \rightarrow \infty}
P_N \left[\ \sup_{0\leq t \leq T} ||M_{t,i}^N||_{-m}>M\right]
=0.$$ We need only
$$\lim_{M \rightarrow \infty}\limsup_{N \rightarrow \infty}
P_N \left[\ ||\int_0^T\tilde Y_{t,i}^N||_{-m}>M\right] =0,$$ again
for $i=1,2$. However, this follows from Proposition
\ref{usefultight3} and Lemma \ref{expbound}.  This completes the
proof of tightness in the colour case.

\bigskip

\noindent\textbf{Acknowledgements.}  This paper is a direct result
of the my dissertation, completed under the supervision of Jeremy
Quastel. I wish to thank him here for suggesting to me this
problem, and for teaching me so much.

\end{document}